\documentclass[10pt]{article}

\usepackage[english]{babel}
\usepackage[a4paper,top=20mm,bottom=20mm,left=25mm,right=25mm,marginparwidth=17.5mm]{geometry}

\usepackage{algorithm,algpseudocode}
\usepackage{amsmath}
\usepackage{amssymb}
\usepackage[backend=biber, style=numeric]{biblatex}
\usepackage{caption}
\usepackage{csquotes}
\usepackage{graphicx}
\usepackage[colorlinks=true, allcolors=blue]{hyperref}
\usepackage{listings}
\usepackage{makecell}
\usepackage[version=4]{mhchem}
\usepackage{microtype}
\usepackage{multirow}
\usepackage{nicematrix}
\usepackage{standalone}
\usepackage{subcaption}
\expandafter\def\csname ver@etex.sty\endcsname{} 
\usepackage{tikz-network}

\usepackage{mycolors}
\usepackage{mycommands}

\usetikzlibrary{positioning}
\title{Automatic partitioning for the low-rank integration of stochastic Boolean reaction networks}

\author{Lukas Einkemmer\thanks{Department of Mathematics, Universit\"{a}t Innsbruck, Innsbruck, Tyrol, Austria (\email{lukas.einkemmer@uibk.ac.at}, \email{julian.mangott@uibk.ac.at})} \and Julian Mangott\footnotemark[1] \and Martina Prugger\thanks{Max-Planck-Institut f\"{u}r Plasmaphysik, Garching, Bavaria, Germany (\email{martina.prugger@ipp.mpg.de})}
}

\begin{document}
\maketitle

\begin{abstract}
Boolean reaction networks are an important tool in biochemistry for studying mechanisms in the biological cell. However, the stochastic formulation of such networks requires the solution of a master equation which inherently suffers from the curse of dimensionality. In the past, the dynamical low-rank (DLR) approximation has been repeatedly used to solve high-dimensional reaction networks by separating the network into smaller partitions. However, the partitioning of these networks was so far only done by hand. In this paper, we present a heuristic, automatic partitioning scheme based on two ingredients: the Kernighan-Lin algorithm and information entropy. Our approach is computationally inexpensive and can be easily incorporated as a preprocessing step into the existing simulation workflow. We test our scheme by partitioning  Boolean reaction networks on a single level and also in a hierarchical fashion with tree tensor networks. The resulting accuracy of the scheme is superior to both partitionings chosen by human experts and those found by simply minimizing the number of reaction pathways between partitions.
\end{abstract}

\section{Introduction}
The biological cell is a tremendously complex system, consisting of thousands of different proteins which can react in many different ways with each other. Gaining insight into these processes is crucial for developing new drugs and methods for cancer treatment. Due to the intricate nature of the reactions, intermediate reaction products are often not accessible by experiment and therefore simulations are needed to complement the experimental observations.

The reactions are often modeled within the framework of chemical reaction systems. The chemical species in such a system can react with each other. Each reaction is characterized by a change in the population number of the affected species and, in the deterministic setting, by a reaction rate. The deterministic approach models the time evolution of the population numbers via rate equations, which amounts to solving a system of ordinary differential equations (ODEs). In the stochastic framework fluctuations, that are inherent in many biological systems, are taken into account. In this approach, the so-called propensity functions replace the reaction rates and instead of rate equations the time evolution of a probability distribution is computed. This probability distribution describes how likely a given number of molecules of a specific chemical species can be found at a given time in the system. The probability distribution in question can be obtained either by solving a master equation or by generating samples of population numbers at various time points according to the stochastic simulation algorithm (SSA) \cite{Gillespie_1976,Gillespie_1977}. The downside of modeling chemical reaction systems, be it with the deterministic or the stochastic description, is that many parameters are required (such as rate constants or propensity functions), for which experimental data are either very noisy or do not exist at all \cite{Barbuti_2020}. 

Therefore it is often more effective to move from the quantitative model of chemical reaction systems to a qualitative description, where population numbers are no longer of interest, but only the activation state of the species is considered. In these so-called Boolean reaction networks, the population numbers are replaced by activation states, which are Boolean variables and describe whether species are activated (i.e., are abundant) or not. The activation state of a species is functionally related to the activation state of other species by logical rules. This approach was first proposed in \cite{Kauffman_1969} and was later formalized in \cite{Thomas_1973}. The biological basis for Boolean reaction networks lies in the fact that many mechanisms in the biological cell exhibit switch-like behavior \cite{Shmulevich_2002}. Consequently, Boolean reaction networks are frequently used to model gene regulatory networks \cite{LeNovere_2015,Barbuti_2020} and have been successfully applied to model the pathways underlying cell fate in cancer \cite{Calzone_2010,Grieco_2013}. Boolean reaction networks have also been useful for inferring regulatory interactions from data about cellular states before and after a perturbation \cite{Prugger_2021}.

The two most common techniques for describing the dynamics of Boolean reaction networks are the synchronous and the asynchronous update schemes, which update the Boolean reaction network in discrete time steps. If the Boolean reaction system is updated synchronously, then in each time step all species are updated at once according to the logical rules of the network.
When using the asynchronous update on the other hand, the species are updated one at a time, and every species is equally likely to be updated, thus resulting in non-deterministic dynamics \cite{Thomas_1991}.

In order to describe the dynamical features of the system, we consider in this paper the stochastic formulation of Boolean reaction networks based on a master equation. This approach is equivalent to the asynchronous update scheme in the limit where time goes to zero \cite{Wang_2012}.

The master equation for Boolean reaction networks cannot be solved analytically for most realistic cases, therefore one has to resort to numerical simulations. However, solving the master equation numerically is a computational problem that suffers from the so-called \emph{curse of dimensionality}. The introduction of each additional species increases the dimensionality of the probability distribution, thus resulting in an exponential growth of the computational time and the memory requirements. This renders a direct simulation of the resulting master equation infeasible for all but very small systems. However, for the description of interesting processes in the biological cell, a large number of species is usually necessary.



A method to overcome the curse of dimensionality is the dynamical low-rank (DLR) approximation. The main idea of this technique is to approximate the high-dimensional probability distribution with a linear combination of a small set of lower-dimensional basis functions, the so-called low-rank factors. Both memory requirements and the computational effort are significantly reduced, since only a small number of basis functions (the so-called \emph{rank of the approximation}) are used. In order to advance the approximation forward in time, the evolution equation for the probability distribution is replaced by evolution equations for the individual degrees of freedom of the low-rank approximation. These evolution equations are obtained by projecting the time derivative of the approximation on the tangent space of the low-rank manifold, which ensures that the approximation remains low-rank \cite{Koch_2007}. Nowadays, the two state-of-the-art methods for the DLR approximation are the projector-splitting \cite{Lubich_2014} and the Basis Update \& Galerkin (BUG) integrators \cite{Ceruti_2022}. They have been improved in the past few years by proposing conservative and asymptotic-preserving schemes \cite{Einkemmer_2024a,Einkemmer_2022,Einkemmer_2021b}, parallel versions \cite{Kusch_2024,Ceruti_2024a} and rank-adaptive \cite{Ceruti_2022a} or second-order \cite{Ceruti_2024b} methods. Both integrators were also extended from matrices to Tucker tensors \cite{Lubich_2018} and tree tensor networks (TTNs) \cite{Ceruti_2020,Ceruti_2023,Ceruti_2024c}. We note that the DLR approximation was used for many applications, such as plasma physics \cite{Einkemmer_2018,Cassini_2021,Coughlin_2022,Einkemmer_2020,Einkemmer_2023,Coughlin_2024}), radiation transport \cite{Peng_2020,Peng_2021,Kusch_2021,Einkemmer_2021,Einkemmer_2021a,Kusch_2022}) and quantum spin systems \cite{Ceruti_2024,Sulz_2024}.


Moreover, the DLR approximation was also successfully applied to chemical reaction systems \cite{Einkemmer_2024,Nicholson_2023,Einkemmer_2024b}, and in particular, to the master equation for Boolean reaction networks \cite{Prugger_2023}. The main idea of the latter paper was to separate the physical dimensions (i.e., the species) of the problem. The reaction network is split into two partitions, and every partition is then described by a small number of low-rank factors. Thus reactions between species in the same partition are treated exactly. An approximation is only performed if a reaction involves species of both partitions, i.e., when the reaction crosses the partition boundary. Ideally, each partition should contain approximately one half of the total species, in order to reduce the memory requirements, but aside from that, there is complete freedom in choosing the partitions. Depending on the concrete problem, species can be kept together to capture important reaction pathways and they can be separated when only a few reactions cross partition boundaries. The DLR approximation is noise-free and resolves rare events more efficiently than Monte Carlo methods such as SSA \cite{Nicholson_2023,Einkemmer_2024}.

Similar approaches for reaction networks based on the DLR approximation have been proposed by \cite{Jahnke_2008,Hegland_2011}, where a low-rank ansatz with basis functions depending only on a single species was used, or by \cite{Kazeev_2014,Kazeev_2015,Dolgov_2015,Dinh_2020,Ion_2021}, where quantized tensor trains (QTTs) have been applied. However, the latter two methods do not have the flexibility of keeping strongly correlated species together. Given the intricate structures of complex biological networks \cite{Barabasi_2009}, it is doubtful that in biological applications we can consider each species independently. This is further supported by the observation that the specific choice of the partitions has a vast impact on the accuracy of the integrator \cite{Prugger_2023}.

In contrast to other equations (for example kinetic equations in plasma physics, where the position coordinates are separated from the velocity coordinates), it is a priori not clear which species should share a partition of the Boolean reaction network. Numerical experiments in \cite{Prugger_2023} indicate that the required rank for a given accuracy of the solution seems to increase with the number of severed reaction pathways and also depends on the importance of these reaction pathways. By choosing an appropriate partitioning, the DLR approximation can therefore be computed with a much smaller rank for a prescribed accuracy than for a ``bad'' partitioning with many cuts. This can drastically reduce the computational and memory requirements. Good partitionings can sometimes be selected by using expert knowledge (as has been done in all papers in the literature). However, for large networks this can be time consuming and error prone. Moreover, it precludes the possibility of rapidly trying out different models as is needed e.g. in model inference problems (as change in the network topology, in general, requires us to change the partitioning as well). Since the number of partitionings scales exponentially in the number of species, exploring all possible partitioning is clearly infeasible. Moreover, we want to avoid running DLR simulations with multiple partitionings as this is potentially expensive.

In this paper, we therefore propose a heuristic scheme for partitioning Boolean reaction networks in an automatized way. Our scheme is based on two ingredients: First, a pool of partitionings is generated with the Kernighan-Lin algorithm \cite{Kernighan_1970}, which is a heuristic graph partitioning algorithm that minimizes the number of severed reaction pathways (edges in the language of graph theory). Second, the information entropy of each partition in the pool is computed and the one with minimal entropy is chosen as the ``best'' partitioning. This algorithm is cheap compared to the time integration of the DLR approximation and thus can be easily integrated into the simulation workflow.

The remainder of this paper is structured as follows. In Section~\ref{sec:boolean-reaction-system} we explain the stochastic formulation of Boolean reaction networks and in Section~\ref{sec:DLRA} give a short overview on the dynamical low-rank approximation. Our heuristic partitioning algorithm based on information entropy and the Kernighan-Lin algorithm is then formulated in Section~\ref{sec:automatic-partitioning}. In Section~\ref{sec:experiments} we investigate the effectiveness for three biochemical examples and compare the automatically found partitionings to the manual choices made in the literature, where possible. We study partitionings where the network is decomposed into two partitions on a single-level and, for the first time in literature, we also perform hierarchical decompositions of Boolean reaction networks. Finally, we conclude in Section~\ref{sec:conclusion}.

\section{Stochastic description of a Boolean reaction network}\label{sec:boolean-reaction-system}
A Boolean reaction network consists of $d$ interacting species $S_0, \dots, S_{d-1}$. In the Boolean setting we only consider two states. That is, a species is either activated (abundant) or deactivated (non-abundant). We call $x_i$ the \emph{activation state} of the $i$-th species $S_i$, so consequently $x_i = 1$ if $S_i$ is activated and $x_i = 0$ if $S_i$ is deactivated. We collect all activation states in the vector $\vb{x} = (x_0,\dots,x_{d-1}) \in \Omega$, with $\Omega = \{0, 1\}^d$.
In the stochastic description the activation states are actually random variables, thus the state of the $i$-th species at time $t$ is given by the random variable $\mathcal{X}_i(t) \in \{0, 1\}$. Therefore, the total activation state $\vb{x}$ is a realization of the random variable $\vb{\mathcal{X}}(t) = (\mathcal{X}_0(t), \dots, \mathcal{X}_{d-1}(t))$ at time $t$.

As mentioned above, the species can react with each other, and the effect of a successful reaction is the change of the activation state of a single species $S_i$. This change of $x_i$ is described by a Boolean rule $\mathcal{B}_i(\vb{x}):\Omega \to \{0, 1\}$. Boolean rules combine the activation states with the logical conjunction ($\land$), the logical disjunction ($\lor$) and the logical negation ($\neg$). After a reaction via the Boolean rule $\mathcal{B}_i(\vb{x})$, the state $\vb{x}$ is modified to $\vb{x}'=(x_0, \dots, x_{i-1}, \mathcal{B}_i(\vb{x}), x_{i+1}, \dots, x_{d-1})$ according to the asynchronous update scheme \cite{Barbuti_2020}. Note that for this scheme each species is activated according to a single Boolean rule, thus, we have in total $d$ different Boolean rules. We also emphasize that the Boolean rule $\mathcal{B}_i(\vb{x})$ modifying $x_i$ is in principle defined on the activation states $\vb{x}$ of all species in the system, but for most reaction systems, $\mathcal{B}_i(\vb{x})$ depends only on a small subset of all species.

Since $\vb{\mathcal{X}}$ is a random variable, we describe the reaction system by the probability density 
\begin{equation*}
    P(t,\vb{x}) = \mathbb{P}(\mathcal{X}_0(t) = x_0, \dots, \mathcal{X}_{d-1}(t) = x_{d-1}),
\end{equation*}
which is the solution of the \emph{master equation for Boolean reaction networks},
\begin{equation}\label{eq:bme}
    \partial_t P(t, \vb{x}) = \sum_{i=0}^{d-1} \sum_{\vb{x}'\in\Omega} \left( A_i(\vb{x}, \vb{x}') P(t, \vb{x}') -  A_i(\vb{x}', \vb{x}) P(t, \vb{x}) \right),
\end{equation}
where $A_i(\vb{x}', \vb{x})$ is the \emph{transition probability per unit time} from state $\vb{x}$ to state $\vb{x}'$ under the Boolean rule $\mathcal{B}_i(\vb{x})$. We require that all reactions (i.e., all Boolean rules) are equally likely and thus we set the transition probability to
\begin{equation*}
    A_i(\vb{x}', \vb{x}) =
    \begin{cases}
        1 & \text{if} \quad \vb{x}' =  (x_0, \dots, x_{i-1}, \mathcal{B}_i(\vb{x}), x_{i+1}, \dots, x_{d-1}),\\
        0 & \text{else}.
    \end{cases}
\end{equation*}
Note that the species in the Boolean reaction networks are assumed to be well-stirred, thus the Boolean rules, the transitions probability and the probability distribution do not show a spatial dependency.

The first term of the master equation is a gain for the probability of state $\vb{x}$ due to transitions from other states $\vb{x}'$ and the second term is a loss due to transitions from $\vb{x}$ into other states \cite{VanKampen_2007}. The master equation can be broken down as a system of ODEs of size $2^d$, with one equation for every state $\vb{x}$.

Boolean networks can be represented as directed graphs. In this representation, each node corresponds to one species and the edge from species $S_j$ to $S_i$ represents the dependency of $S_i$ on $S_j$ due to the Boolean rule $\mathcal{B}_i(\vb{x})$. We call these edges also \emph{reaction pathways}. Since each Boolean rule affects only the activation of a single species, there can be at maximum two edges with a different direction between two nodes.

\begin{figure}[!htb]
    \begin{minipage}[b]{0.47\textwidth}
        \centering
        \begin{tabular}{c|c}
            \hline
            $i$ & Boolean rule $\mathcal{B}_i(\vb{x})$
            \tabularnewline
            \hline 
            $0$ & $\neg x_1$
            \tabularnewline
            $1$ & $x_2 \lor x_3$
            \tabularnewline
            $2$ & $x_2 \lor (x_2 \land x_4)$
            \tabularnewline
            $3$ & $(\neg x_2) \land x_3$
            \tabularnewline
            $4$ & $x_1 \land (\neg x_4)$
            \tabularnewline
            \hline 
        \end{tabular}
        \captionof{table}{An example for a Boolean reaction network with five species. The corresponding graph is shown in Figure~\ref{fig:example-graph}.}
        \label{tab:example-rules}
    \end{minipage}
    \hfill
    \begin{minipage}[b]{0.47\textwidth}
        \centering
        \includegraphics[height=0.55\textwidth]{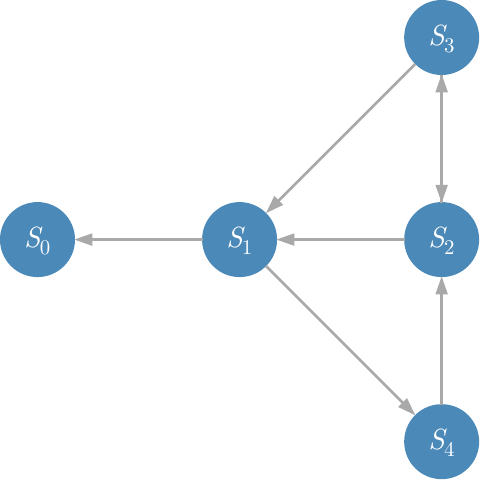}
        \captionof{figure}{The graph of the Boolean reaction network with five species and rules specified in Table~\ref{tab:example-rules}.\label{fig:example-graph}}
    \end{minipage}
\end{figure}

We want to elucidate the concepts of Boolean reaction networks by giving a small example. Let us consider a reaction network with five species $S_i$ ($i=0,\dots,4$) and Boolean rules as stated in Table~\ref{tab:example-rules}. Then the corresponding graph is shown in Figure~\ref{fig:example-graph}. Note that in principle we would also have to draw self-loops for $S_2$, $S_3$ and $S_4$. However, it will become clear in Section~\ref{sec:automatic-partitioning} that we are interested in partitions of the reaction network which minimize the number of connections between the nodes; therefore, self-loops are not relevant for our purposes and hence we do not consider them in the reaction graphs. We also want to emphasize that the graph of a Boolean reaction network only contains information about the dependencies between species, but it does not yield any information about the concrete form of the Boolean rules. Therefore, two different Boolean reaction networks could share the same graph representation.

We proceed with the example by calculating the transition probability for a single rule. The transition probability due to the Boolean rule $\mathcal{B}_0(\vb{x})$ can be obtained by regarding the corresponding truth table shown in Table~\ref{tab:rule-0}: We first note that this Boolean rule only depends on species $S_0$ and $S_1$. In the truth table a change of $x_0$ due to $\mathcal{B}_0(\vb{x})$ is only given in the first and the last line -- this means that there are only two possible transitions between states: First, species $S_0$ can be activated if $S_1$ is deactivated, which corresponds to a transition from state $\vb{x}=(0, 0, \cdot, \cdot, \cdot)$ to $\vb{x}'=(1, 0, \cdot, \cdot, \cdot)$. Here, ``$\cdot$" denotes that all other species can be in any activation state and they remain in this state under the transition, since $\mathcal{B}_0(\vb{x})$ does not depend on them. Second, species $S_0$ can be deactivated when $x_1=1$ is activated, which amounts to a transition from state $\vb{x}=(1, 1, \cdot, \cdot, \cdot)$ to $\vb{x}'=(0, 1, \cdot, \cdot, \cdot)$.

\begin{table}[!htb]
    \begin{subtable}[t]{0.3\textwidth}
        \centering
        \begin{tabular}{c|c|c}
            \hline
            $x_0$ & $x_1$ & $\mathcal{B}_0(\vb{x})$
            \tabularnewline
            \hline 
            0 & 0 & 1
            \tabularnewline
            0 & 1 & 0
            \tabularnewline
            1 & 0 & 1
            \tabularnewline
            1 & 1 & 0
            \tabularnewline
            \hline
        \end{tabular}
        \caption{\label{tab:rule-0}}
    \end{subtable}
    \hfill
    \begin{subtable}[t]{0.3\textwidth}
        \centering
        \begin{tabular}{c|c|c}
            \hline
            $x_2$ & $x_4$ & $\mathcal{B}_2(\vb{x})$
            \tabularnewline
            \hline 
            0 & 0 & 0
            \tabularnewline
            0 & 1 & 0
            \tabularnewline
            1 & 0 & 1
            \tabularnewline
            1 & 1 & 1
            \tabularnewline
            \hline
        \end{tabular}
        \caption{\label{tab:rule-2}}
    \end{subtable}
    \hfill
    \begin{subtable}[t]{0.3\textwidth}
        \centering
        \begin{tabular}{c|c|c}
            \hline
            $x_2$ & $x_3$ & $\mathcal{B}_3(\vb{x})$
            \tabularnewline
            \hline 
            0 & 0 & 0
            \tabularnewline
            0 & 1 & 1
            \tabularnewline
            1 & 0 & 0
            \tabularnewline
            1 & 1 & 0
            \tabularnewline
            \hline
        \end{tabular}
        \caption{\label{tab:rule-3}}
    \end{subtable}
    \caption{Truth tables for $\mathcal{B}_0(\vb{x})$, $\mathcal{B}_2(\vb{x})$ and  $\mathcal{B}_3(\vb{x})$ of the five-dimensional Boolean reaction network example of Table~\ref{tab:example-rules}.}
\end{table}

The transition probability $A_0(\vb{x}', \vb{x})$ is therefore represented in matrix form as
\begin{equation*}
    A_0 = \ \;
    \begin{NiceArray}{w{c}{0.1cm}>{\scriptstyle \color{gray}}ccccc}[baseline=line-5]
        & & \Block{1-4}{\scriptstyle \color{darkgray} \centering (x_0,x_1)}\\
        \RowStyle{\scriptstyle \color{gray}} 
        & & (0,0) & (0,1) & (1,0) & (1,1)\\
        \Block{4-1}{\scriptstyle \color{darkgray} \rotate{(x_0',x_1')}} & (0,0) & 0 & 0 & 0 & 0\\
        & (0,1) & 0 & 0 & 0 & 1\\
        & (1,0) & 1 & 0 & 0 & 0\\
        & (1,1) & 0 & 0 & 0 & 0
        \CodeAfter
        \SubMatrix[{3-3}{6-6}]
    \end{NiceArray}
    \ \,\otimes
    I_2 \otimes I_3 \otimes I_4,
\end{equation*}
where $I_i$ is the $2 \times 2$~identity matrix in the subspace of species $S_i$.
In a similar fashion, the transition probabilities for the other Boolean rules can be obtained and subsequently, Equation~(\ref{eq:bme}) can be solved. 

For more details about the stochastic formulation of reaction networks and master equations we refer the reader to \cite{VanKampen_2007,Gardiner_2004,Erban_2020}. The reader might also consult \cite{Barabasi_2009,LeNovere_2015} for further information on Boolean reaction networks.

\section{Dynamical low-rank approximation}\label{sec:DLRA}
We mentioned already in the introduction that the master equation and its solution suffer from the curse of dimensionality: The memory requirements scale exponentially with the number of species. Therefore, we have to reduce the system size in order to solve Equation~(\ref{eq:bme}) on currently available hardware. But memory requirements will also be a pressing issue in the future, due to the domain scientists' need for simulations of ever larger and more complex systems.
Moreover, for the time evolution of the full probability distribution with the Boolean master equation, the right-hand side of the equation has to be evaluated and all degrees of freedom have to be updated in each time step. By reducing the number of degrees of freedom, both memory and computational time can be saved.

In \cite{Prugger_2023,Einkemmer_2024,Einkemmer_2024b} this was done by using a dynamical low-rank (DLR) approximation. In the following, we will give a brief overview of the DLR approximation. The main idea of this technique is to separate the species (and thus the reaction network and its corresponding reaction graph) into partitions. This can either be done on only one level, or more generally, on multiple levels in a hierarchical fashion. We first describe the single (``one-level") partitioning of the reaction network.

\subsection{One-level partitioning}
Let $\mathcal{P}=\{S_0,\dots,S_{d-1}\}$ be the set of all species, with $|\mathcal{P}|=d$. Then we perform a single (``one-level") separation of the reaction network by decomposing $\mathcal{P}$ into two smaller partitions $\mathcal{P}^0 \subset \mathcal{P}$ and $\mathcal{P}^1 \subset \mathcal{P}$, such that $\mathcal{P} = \mathcal{P}^0 \cup \mathcal{P}^1$ and $\mathcal{P}^0 \cap \mathcal{P}^1 = \emptyset$. We formally write such a partitioning as $\mathcal{P} = (\mathcal{P}^0, \mathcal{P}^1)$. Analogously, we write a single state as $\vb{x} = (\vb{x}^0, \vb{x}^1)$, where $\vb{x}^0 \in \Omega^0 = \{0,1\}^{d^0}$ and $\vb{x}^1 \in \Omega^1 = \{0,1\}^{d^1}$ are the population numbers of the species in $\mathcal{P}^0$ and $\mathcal{P}^1$, respectively, and $d^0 + d^1 = d$. Therefore, the size of the state space is $|\Omega^0|=2^{d^0}$ for $\mathcal{P}^0$ and $|\Omega^1|=2^{d^1}$ for $\mathcal{P}^1$. In the following, we will denote all quantities belonging to a specific partition (here, $0$ or $1$) with superscripts.

The probability distribution $P(t,\vb{x})$, which is the solution of the master equation for Boolean reaction networks, is then approximated by $r$ time-dependent low-rank factors $X^0_{i_0}(t, \vb{x}^0)$ and $X^1_{i_1}(t, \vb{x}^1)$ and a time-dependent, invertible coefficient matrix $S_{i_0i_1}(t) \in \mathbb{R}$:
\begin{equation}\label{eq:DLR-approximation}
    P(t, \vb{x}) \approx \sum_{i_0,i_1=0}^{r-1} X^0_{i_0}(t, \vb{x}^0) S_{i_0i_1}(t) X^1_{i_1}(t, \vb{x}^1).
\end{equation}
Note that the low-rank factors can be represented as matrices without any further loss of accuracy, since the activation states $\vb{x}^0$ and $\vb{x}^1$ are discrete. Therefore, Equation~(\ref{eq:DLR-approximation}) is similar to a truncated singular value decomposition (SVD), with the main difference that the coefficient matrix $S_{i_0i_1}(t)$ must not necessarily be diagonal.

The number of low-rank factors $r$ is the so-called rank. Recall that the total number of degrees of freedom for the full model is $2^d$. Since the $r$ low-rank factors $X^0_i(t,\vb{x}^0)$ and $X^1_j(t,\vb{x})$ only depend on the species in their partition, the DLR approximation needs in total $r(r + 2^{\frac{d}{2}+1})$ degrees of freedom, when $d^0=d^1=\frac{d}{2}$. If the number of low-rank factors needed to approximate $P(t, \vb{x})$ (up to a given precision) is small compared to the partition sizes, $r \ll |\Omega^0|$ and $r \ll |\Omega^1|$, then the DLR approximation yields a drastic reduction of the memory requirements from $\mathcal{O}(2^d)$ to $\mathcal{O}(2^{\frac{d}{2}})$. Since Boolean rules typically depend only on a small subset of all species, the transition probabilities $A_i(\vb{x}',\vb{x})$ and therefore also the master equation for Boolean reaction networks have a low-rank structure. In Section~\ref{sec:experiments} we will see indeed that a rank of approximately $r=20$ is sufficient for many large Boolean reaction networks.

\begin{figure}[!htb]
    \centering
    \begin{minipage}[t]{0.47\textwidth}
        \begin{subfigure}[b]{\textwidth}
            \centering
            \scalebox{0.95}{\begin{tikzpicture}
    \SetVertexStyle[TextFont=\normalsize]
    \SetEdgeStyle[LineWidth=1.0pt]
    \Vertex[label={$S_{i_0 i_0}$},position=above,style={color=custom_red}]{root}
    \Vertex[x=-1.05,y=-1,label={$X^0_{i_0}(x_0,x_1)$},position=below,style={color=custom_blue}]{0}
    \Vertex[x=1.05,y=-1,label={$X^1_{i_1}(x_2,x_3,x_4)$},position=below,style={color=custom_violet}]{1}
    
    \Edge(root)(0)
    \Edge(root)(1)
\end{tikzpicture}}
            \caption{\label{fig:example-tree-p0}}
        \end{subfigure}
        \\[3ex]
        \begin{subfigure}[b]{\textwidth}
            \centering
            \includegraphics[height=0.75\textwidth]{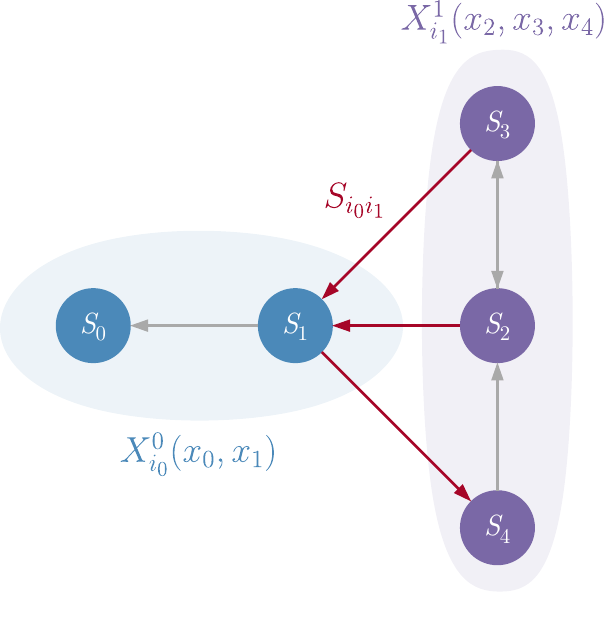}
            \caption{\label{fig:example-graph-p0}}
        \end{subfigure}
        \caption{(a) Graphical representation of Equation~(\ref{eq:DLR-approximation}) and (b) reaction graph for the five-dimensional example problem of Section~\ref{sec:boolean-reaction-system} with two partitions. The time dependency of the low-rank factors and the coefficient matrix has been omitted.}
    \end{minipage}
    \hfill
    \begin{minipage}[t]{0.47\textwidth}
        \begin{subfigure}[b]{\textwidth}
            \centering
            \scalebox{0.95}{\begin{tikzpicture}
    \SetVertexStyle[TextFont=\normalsize]
    \SetEdgeStyle[LineWidth=1.0pt]
    \Vertex[label={$S_{i_0 i_1}$},position=above,style={color=custom_red}]{root}
    \Vertex[x=-2.0,y=-1,label={$Q^0_{i_0 i_{00} i_{01}}$},position=left,style={color=custom_blue}]{0}
    \Vertex[x=2.0,y=-1,label={$Q^1_{i_1 i_{10} i_{11}}$},position=right,style={color=custom_violet}]{1}
    \Vertex[x=-3.0,y=-2,label={$X^{00}_{i_{00}}(x_0)$},position=below,style={color=custom_darkblue}]{00}
    \Vertex[x=-1.0,y=-2,label={$X^{01}_{i_{01}}(x_1)$},position=below,style={color=custom_lightblue}]{01}
    \Vertex[x=1.0,y=-2,label={$X^{10}_{i_{10}}(x_2,x_3)$},position=below,style={color=custom_darkviolet}]{10}
    \Vertex[x=3.0,y=-2,label={$X^{11}_{i_{11}}(x_4)$},position=below,style={color=custom_lightviolet}]{11}
    
    \Edge(root)(0)
    \Edge(root)(1)
    \Edge(0)(00)
    \Edge(0)(01)
    \Edge(1)(10)
    \Edge(1)(11)
\end{tikzpicture}}
            \caption{\label{fig:example-tree-p1}}
        \end{subfigure}
        \\[3ex]
        \begin{subfigure}[b]{\textwidth}
            \centering
            \includegraphics[height=0.75\textwidth]{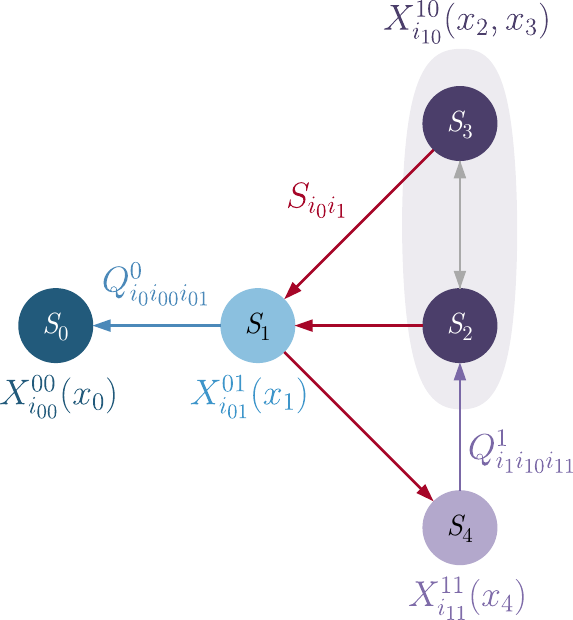}
            \caption{\label{fig:example-graph-p1}}
        \end{subfigure}
        \caption{(a) Graphical representation of Equation~(\ref{eq:DLR-approximation}), with low-rank factors being decomposed according to Equation~(\ref{eq:TTN-X}) and (b) reaction graph for the example problem of Section~\ref{sec:boolean-reaction-system} with two subpartitions for each of the two partitions. For sake of simplicity, the time dependency of the low-rank factors, the coefficient matrix and the connection tensors have been omitted.}
    \end{minipage}
\end{figure}

In Figure~\ref{fig:example-graph-p0} we schematically depict the separation of the five-dimensional example system introduced in Section~\ref{sec:boolean-reaction-system} into two partitions: The partitioning was chosen as $\mathcal{P}^0=\{S_0, S_1\}$ (blue) and $\mathcal{P}^1=\{S_2, S_3, S_4\}$ (purple). Since the low-rank factors $X^0_{i_0}(t,x_0,x_1)$ only depend on the species in partition $\mathcal{P}^0$, it is ensured that reactions depending only on the species in this partition are treated exactly (gray arrows). Similarly, reactions depending only on species which lie in partition $\mathcal{P}^1$ are fully described by the low-rank factors $X^1_{i_1}(t,x_2,x_3,x_4)$. The coefficient matrix $S_{i_0i_1}(t)$ couples the low-rank factors of the different partitions, and therefore approximates the reaction pathways which cross the two partitions (red arrows). Ideally, the partitions are chosen such that the number of approximated reaction pathways is minimized and tightly coupled species are kept together in a single partition. We will see later that this heavily influences the required rank for a given accuracy.

Equation~(\ref{eq:DLR-approximation}) can also be represented in a graphical way. The tree structure of Figure~\ref{fig:example-tree-p0} is obtained, when the low-rank factors and the coefficient matrix are depicted as nodes which are connected according to the contraction of the indices in Equation~(\ref{eq:DLR-approximation}). The coefficient matrix sits on the root and the two low-rank factors are represented by the leaves of this binary tree with height $1$. In the next subsection, we will extend this structure to binary trees of an arbitrary height.

\subsection{Hierarchical partitioning}
If the computation of a Boolean reaction network with a given number of species is barely doable by solving the master equation directly, then the one-level partitioning essentially enables to simulate systems with twice the number of species. This means in turn that for large networks, the state space of even a single partition may become prohibitively large. Therefore, we want to separate partitions into smaller ones, until the sizes of the resulting partitions are computationally tractable.

We illustrate this by subdividing $\mathcal{P}^0$ into two subpartitions $\mathcal{P}^{00} \subset \mathcal{P}^0$ and $\mathcal{P}^{01} \subset \mathcal{P}^0$ with $\mathcal{P}^{00} \cap \mathcal{P}^{01} = \emptyset$, by writing formally $\mathcal{P}^0 = (\mathcal{P}^{00}, \mathcal{P}^{01})$. This amounts to a decomposition of the low-rank factor $X^0_i(t,\vb{x}^0)$ into a connection tensor $Q^0_{ii_{00}i_{01}}$ and two low-rank factors $X^{00}_{i_{00}}(t, \vb{x}^{00})$ and $X^{01}_{i_{01}}(t, \vb{x}^{01})$, where $\vb{x}^0 = (\vb{x}^{00}, \vb{x}^{01})$:
\begin{equation}\label{eq:TTN-X0}
    X^0_{i_0}(t,\vb{x}^0) = \sum_{i_{00}=0}^{r^{00}-1} \sum_{i_{01}=0}^{r^{01}-1} X^{00}_{i_{00}}(t,\vb{x}^{00}) Q^0_{i_0 i_{00} i_{01}}(t) X^{01}_{i_{01}}(t,\vb{x}^{01}).
\end{equation}
This decomposition cannot only be applied to the low-rank factor $X^0_{i_0}(t, \vb{x}^0)$ (and, slightly modified, also to $X^1_{i_1}(t, \vb{x}^1)$), but it can also be adapted for the low-rank factors of the new subpartitions $\mathcal{P}^{00}$ and $\mathcal{P}^{01}$. For an arbitrary partition $\tau$ with $\vb{x}^\tau = (\vb{x}^{\tau_0}, \vb{x}^{\tau_1})$ and subpartitions $\tau_0$ and $\tau_1$, Equation~(\ref{eq:TTN-X0}) can be generalized to
\begin{equation}\label{eq:TTN-X}
    X^\tau_{i_{\tau}}(t,\vb{x}^{\tau}) = \sum_{i_{\tau_0}=0}^{r^{\tau_0}-1} \sum_{i_{\tau_1}=0}^{r^{\tau_1}-1} X^{\tau_0}_{i_{\tau_0}}(t,\vb{x}^{\tau_0}) Q^\tau_{i_\tau i_{\tau_0} i_{\tau_1}}(t) X^{\tau_1}_{i_{\tau_1}}(t,\vb{x}^{\tau_1}).
\end{equation}
In contrast to the SVD-like decomposition of the one-level partitioning, Equations~(\ref{eq:TTN-X0}) and (\ref{eq:TTN-X}) are so-called Tucker decompositions and $X^\tau_{i_\tau}(t,\vb{x}^\tau) \equiv X^\tau_{i_\tau}(t,\vb{x}^{\tau_0},\vb{x}^{\tau_1})$ is treated as a time-dependent $r^\tau \times |\Omega^{\tau_0}| \times |\Omega^{\tau_1}|$~tensor. Remember that for the matrix case the number of low-rank factors $X^0_{i_0}(t,\vb{x}^0)$ and $X^1_{i_1}(t,\vb{x}^1)$ was the same, namely $r$. This spared us from introducing different ranks $r^0$ and $r^1$ for the low-rank factors $X^0_{i_0}(t,\vb{x}^0)$ and $X^1_{i_1}(t,\vb{x}^1)$, respectively. For the Tucker decomposition the situation is different: Since the connection tensor $Q^\tau_{i_\tau i_{\tau_0} i_{\tau_1}}(t)$ with full multilinear rank $(r^\tau, r^{\tau_0}, r^{\tau_1})$ has to obey the condition that \cite{Ceruti_2020}
\begin{equation}\label{eq:Tucker-condition}
    (r^\tau \le r^{\tau_0} r^{\tau_1}) \land (r^{\tau_0} \le r^{\tau_1} r^\tau) \land (r^{\tau_1} \le r^{\tau_0} r^\tau),
\end{equation}
the ranks $r^{\tau_0}$ and $r^{\tau_1}$ can be different from each other and therefore also the number of low-rank factors $X^{\tau_0}_{i_{\tau_0}}(t,\vb{x}^{\tau_0})$ and $X^{\tau_1}_{i_{\tau_1}}(t,\vb{x}^{\tau_1})$.

With Equation~(\ref{eq:TTN-X}), the system can be decomposed hierarchically until the system size is reduced sufficiently. Theoretically, this decomposition could be done until the low-rank factors only depend on a single species. In our view it is however often preferable to keep species which are coupled by biologically important reaction pathways together in a single partition.

In Figure~\ref{fig:example-graph-p1} we demonstrate the hierarchical partitioning for our five-dimensional example problem from Section~\ref{sec:boolean-reaction-system}. Each of the two partitions is now subdivided into two smaller subpartitions. The species $S_0$ and $S_1$ are lying in different subpartitions and the reaction between these species is approximated by the connection tensor $Q^0_{i_0 i_{00} i_{01}}(t)$ (blue arrows). Reactions between the large partitions are again approximated via the coefficient matrix $S_{i_0 i_1}(t)$ (red arrows), and the pathway between species $S_4$ and $S_2$ is approximated by the connection tensor $Q^1_{i_1 i_{10} i_{11}}(t)$ (purple arrow). Reactions between species $S_2$ and $S_3$ are treated exactly (gray arrows) by the low-rank factor $X^{10}_{i_{10}}(t, x_2, x_3)$.

The hierarchical partitioning of this five-dimensional example problem with Equations~(\ref{eq:DLR-approximation}) and (\ref{eq:TTN-X}) can be represented as a binary tree of height $2$ (Figure~\ref{fig:example-tree-p1}). As for the one-level partitioning, the coefficient matrix again sits on the root and the low-rank factors on the leaves of the binary tree, but now the tree has also internal nodes, which are associated with the connection tensors.

In general, any hierarchical partitioning using Equations~(\ref{eq:DLR-approximation}) and (\ref{eq:TTN-X}) can be represented as a binary tree, and each (sub-)partition $\tau$ can be identified with a node of the binary tree. The resulting decomposition is similar to a higher-order SVD (HOSVD) and the structure of the decomposition is also known as a \emph{tree tensor network} (TTN). If $\tau$ is not further subdivided into smaller partitions, then it is a leaf node of the tree, but if $\tau$ has subpartitions $\tau_0$ and $\tau_1$, then $\tau$ corresponds to an internal node with child nodes $\tau_0$ and $\tau_1$. Leaf nodes are associated with low-rank factors, internal nodes with connection tensors.

\subsection{Time integration}
In the previous section we have seen that a partitioning yields in general a TTN, which can be depicted as a binary tree. The root of the tree corresponds to a coefficient matrix, the internal nodes to connection tensors and the leaves to low-rank factors. The probability distribution can be reconstructed for any time $t$ from these quantities via Equations~(\ref{eq:DLR-approximation}) and (\ref{eq:TTN-X}).
Yet, forming the entire probability distribution $P(t,\vb{x})$ for the time evolution is not a viable option: In order to overcome the curse of dimensionality, we have to replace the master equation with individual evolution equations for the coefficient matrix, the connection tensors and the low-rank factors. These evolution equations have been derived in \cite{Prugger_2023,Einkemmer_2024} for one-level partitioning by using the projector-splitting integrator proposed in \cite{Lubich_2014}. Furthermore, the evolution equations have been generalized for hierarchical partitioning in \cite{Einkemmer_2024b} by using an extension of the projector-splitting integrator to tree tensor networks proposed in \cite{Ceruti_2020}. We will not go into more details here as it is not relevant for the automatic partitioning algorithm. The interested reader is referred to the listed publications in which details are provided.

For the time integration the initial condition has to be approximated by
\begin{equation*}
    P(0,\vb{x}) \approx \sum_{i_0,i_1=0}^{r-1} X^0_{0,i_0}(\vb{x}^0) S_{0,i_0 i_1} X^1_{0,i_1}(\vb{x}^1),
\end{equation*}
and, in the case of a hierarchical partitioning, the low-rank factors have to be recursively approximated by
\begin{equation*}
    X^\tau_{0,i_\tau}(\vb{x}^\tau) \approx \sum_{i_{\tau_0}=0}^{r^{\tau_0}-1} \sum_{i_{\tau_1}=0}^{r^{\tau_1}-1} X^{\tau_0}_{0,i_{\tau_0}}(\vb{x}^{\tau_0}) Q^{\tau}_{0,i_\tau i_{\tau_0} i_{\tau_1}} X^{\tau_1}_{0,i_{\tau_1}}(\vb{x}^{\tau_1}).
\end{equation*}
When $P(t, \vb{x})$ is already in the low-rank form, then obtaining the coefficient matrix, the connection tensors and the low-rank factors is straightforward. Otherwise these quantities can be computed from $P(0, \vb{x})$ by a truncated randomized (HO)SVD.

\section{Automatic partitioning}\label{sec:automatic-partitioning}
In Section~\ref{sec:DLRA} we mentioned that the required rank for a prescribed accuracy of the DLR approximation depends on the amount and the concrete type of the approximated reaction pathways. Therefore, by using a ``good" partitioning, we can reduce the required rank for a given accuracy and thereby memory and computational time can be reduced.

In the literature partitionings have so far been chosen only manually  \cite{Prugger_2023,Einkemmer_2024,Einkemmer_2024b}. This has been guided by minimizing cuts of reaction pathways and usually also involves some expert knowledge (e.g.~to identify tightly coupled parts of the network that should be kept within a partition). For larger networks this can be a challenging and tedious task. Moreover, it is often unclear whether a given partitioning works well. For example, in \cite{Prugger_2023} a number of ways to partition the network have been compared and it is observed that even carefully chosen partitionings can have relatively large differences in accuracy (for a fixed rank). Using a hierarchical decomposition further complicates the task of finding a good partitioning.

Thus, it would be desirable to automate this process. In addition to freeing the user from manually finding a good partitioning, it also allows the use of low-rank methods in cases where a number of different reaction networks (i.e.~models) are explored, for example with an optimization algorithm as is commonly done in model inference problems. Ideally, this automatic scheme should be computationally cheap compared to the actual time integration of the DLR approximation. Clearly, it is not feasible to explore all possible network configurations (the number of network scales exponentially in the number of species \cite{Kernighan_1970}) or to run the low-rank simulation with a number of selected partitionings to determine their utility. What we do instead is to propose a heuristic that evaluates partitionings by a) the number of cuts in reaction pathways made and b) how important these cuts are. The latter is quantified by the loss of information due to the cut. Before this is done, a number of suitable candidates with a small number of cuts (a ``partitioning pool") is identified by the Kernighan-Lin algorithm.


In the next two sections, we will present the two main ingredients of our partitioning algorithm for the special case of one-level partitionings; later in Section~\ref{sec:algorithm} we will extend our method to general hierarchical partitionings.


\subsection{Kernighan-Lin algorithm}
Our first strategy of separating the Boolean reaction system into two parts is by minimizing the severed reaction pathways of the corresponding reaction graph. The Kernighan-Lin algorithm \cite{Kernighan_1970} does exactly that: This heuristic algorithm separates an undirected graph with numerical weights on the edges into two partitions of approximately equal size.

Recall that the graph of a Boolean reaction network is a directed graph. So before performing the Kernighan-Lin algorithm, we first have to convert the reaction graph into an undirected graph by taking the number of reaction pathways between two species as weights. Since for Boolean reaction networks two nodes are connected by at most two edges of different direction, a weight can be either $1$ or $2$. At the moment we are only interested in minimizing the number of reaction pathways between partitions, therefore we can omit the direction of the reaction pathways.

The Kernighan-Lin algorithm starts with an arbitrary partitioning $(\mathcal{P}^0, \mathcal{P}^1)$ of the reaction network with approximately equal size, with $d^0 = \lceil \frac{d}{2}\rceil$ species in $\mathcal{P}^0$ and $d^1 = \lfloor \frac{d}{2} \rfloor$ in $\mathcal{P}^1$. Thereby, the total number of degrees of freedom (and thus the memory requirements) for the DLR approximation are minimized.


Next, the Kernighan-Lin algorithm tries to decrease the cost (which is the sum of all weighted edges crossing the two partitions) of the partitioning by subsequently swapping pairs of nodes $(S_i,S_j)$, where $S_i \in \mathcal{P}^0$ and $S_j \in \mathcal{P}^1$. A local minimum is found when no further improvement is possible by interchanging $S_i$ and $S_j$. 

Before swapping a pair of nodes, the Kernighan-Lin algorithm, in its standard implementation, requires that the nodes are sorted according to a cost function. This sorting step has typically complexity $\mathcal{O}(d\log{d})$. Since all $d$ nodes of the graph are swapped by the algorithm, the overall computational cost scales non-exponentially with $\mathcal{O}(d^2\log{d})$.


Since the Kernighan-Lin algorithm is heuristic and therefore does not guarantee to find the global minimum, we generate possible candidates by running the algorithm multiple times with random initial partitionings. Due to the favorable complexity, the overhead compared to the time integration of the DLR approximation is negligible even for large networks (such as the $41$-dimensional apoptosis network in Section~\ref{sec:experiments}).

The main problem with this approach is that only the topology (i.e., the graph) of the reaction network is taken into account. In Section~\ref{sec:boolean-reaction-system}, we mentioned that in principle two different reaction networks could exhibit the same reaction graph, and in this case the Kernighan-Lin algorithm would yield the same partitioning results for both networks. However, it was already demonstrated in \cite{Prugger_2023} that minimizing the reaction pathways between partitions does not always lead to the ``best" partitioning, i.e., the one which needs the lowest rank for a given accuracy of the DLR approximation. The Kernighan-Lin algorithm may be a good starting point, but in order to find better partitionings, we also have to consider information about the concrete type of the Boolean rules. We will address this issue in the following section.

\subsection{Information entropy}\label{sec:information-entropy}
The second ingredient of our automatic partitioning scheme will address the shortcomings of the Kernighan-Lin algorithm by quantifying the information loss of the approximated reaction pathways. Although the transition probability has the same form for all Boolean rules, they are not all equally important. We want to elucidate this by revisiting the first three Boolean rules of the five species example considered in Section~\ref{sec:boolean-reaction-system}.

Let us first regard the Boolean rule $\mathcal{B}_0(\vb{x})$, whose truth table is shown in Table~\ref{tab:rule-0}. The rule activates $S_0$ if $S_1$ is deactivated, and vice versa. Therefore, the activation state of $S_0$ depends completely on $S_1$. If the pathway between $S_0$ and $S_1$ is cut, then all the information for the change of the activation state of $S_0$ is lost and thus the required rank increases.

Next, we consider the Boolean rule $\mathcal{B}_2(\vb{x})$, which in mathematical logic is also known as the ``absorption law". When looking at the truth table shown in Table~\ref{tab:rule-2}, the name of this rule becomes immediately clear: The theoretical dependency of $S_4$ is completely absorbed and thus the activation of $S_2$ actually does not depend on $S_4$. Cutting this pathway should not increase the rank.

Finally, we study the Boolean rule $\mathcal{B}_3(\vb{x})$, whose truth table is shown in Table~\ref{tab:rule-3}. If species $S_3$ is deactivated, then it always stays deactivated, regardless of the activation state of $S_2$. Only when $S_3$ is activated, then $S_2$ has an influence on the activation state of $S_3$. Therefore, not all information is lost when cutting the reaction pathway, and the information loss should lie between the one observed for $\mathcal{B}_0(\vb{x})$ (complete loss) and $\mathcal{B}_2(\vb{x})$ (no loss).

We quantify the information loss of the approximated pathways by the notion of information entropy, which we will define in the following. Let us consider an arbitrary partitioning $(\mathcal{P}^0, \mathcal{P}^1)$ and a species $S_i$ whose activation state is changed by the Boolean rule $\mathcal{B}_i(\vb{x})$. Next, we identify the partition which contains species $S_i$ with $\mathcal{P}_i$ and the other partition with $\mathcal{P}^\complement_i$, so a state $\vb{x}$ can also be written as $\vb{x} = (\vb{x}_i, \vb{x}_i^\complement)$, with $\vb{x}_i \in \Omega_i$ and $\vb{x}_i^\complement \in \Omega_i^\complement$. For example, we have $\mathcal{P}_i=\mathcal{P}^0$ and $\mathcal{P}_i^\complement = \mathcal{P}^1$ if $S_i \in \mathcal{P}^0$. We emphasize that this notation is necessary, because (in contrast to the Kernighan-Lin algorithm) the information entropy depends on the direction of the reaction pathways (information is lost from the viewpoint of the species whose activation state is changed).

Then, the loss of information by cutting a reaction pathway related to $\mathcal{B}_i(\vb{x})$ is measured heuristically with the rule-specific entropy
\begin{equation}\label{eq:entropy-rule}
    h_i = \frac{1}{|\Omega_i|}\sum_{\vb{x}_i \in \Omega_i} \tilde{h}_i(\vb{x}_i),
\end{equation}
where
\begin{equation}\label{eq:entropy-state}
    \tilde{h}_i(\vb{x}_i)= -p_i(\vb{x}_i) \log_2 p_i(\vb{x}_i) - (1-p_i(\vb{x}_i)) \log_2(1-p_i(\vb{x}_i))
\end{equation}
is the information entropy introduced by Shannon \cite{Shannon_1948} and $p_i(\vb{x}_i)$ the probability that $\mathcal{B}_i(\vb{x}_i,\vb{x}_i^\complement)=1$ for a given $\vb{x}_i$, i.e.
\begin{equation}\label{eq:probability}
    p_i(\vb{x}_i) = \frac{1}{|\Omega^\complement_i|} \sum_{\vb{x}_i^\complement \in \Omega_i^\complement} \mathcal{B}_i(\vb{x}_i,\vb{x}_i^\complement).
\end{equation}
The total entropy $H$ of the reaction network (and thus the total loss of information) for the partitioning $(\mathcal{P}^0,\mathcal{P}^1)$ is given by
\begin{equation}\label{eq:entropy}
    H = \sum_{i=0}^{d-1} h_i.
\end{equation}
Note that in contrast to $h_i(\vb{x}_i)$, the total entropy $H$ is not normalized, and therefore can assume values in the interval $[0, d]$.

Let us now quantify the loss of information for the three exemplary rules from above: The rule-specific entropy for $\mathcal{B}_0(\vb{x})$ is $h_0=1$, so we indeed loose the full information when cutting the associated reaction pathways. For $\mathcal{B}_2(\vb{x})$ the entropy is $h_2=0$, so no information is lost. Finally, for $\mathcal{B}_3(\vb{x})$, we obtain $h_3=\frac{1}{2}$ and therefore half of the information is lost. In total, the results for the rule-specific entropy exactly match our expectations from above.

From a computational perspective, the most expensive part is the calculation of the probabilities $p_i(\vb{x}_i)$: Let $d_i^\mathrm{dep}$ be the number of species on which the Boolean rule $\mathcal{B}_i(\vb{x})$ actually depends. In Section~\ref{sec:boolean-reaction-system} we noted that usually the Boolean rules only depend on a small subset of all of the species, so $d_i^\mathrm{dep} \ll d$. For instance,  $d_i^\mathrm{dep}$ is at most $6$ for all of the example problems presented in Section~\ref{sec:experiments}. Then the computational complexity for calculating $p_i(\vb{x}_i)$ is $\mathcal{O}(2^{d_i^\mathrm{dep}})$.

\subsection{Algorithm}\label{sec:algorithm}
In this section, we will put the concepts of the previous two sections together and propose a  partitioning scheme for Boolean reaction networks.

We first consider the special case of the one-level partitioning. The proposed algorithm has essentially two steps: First, we generate a pool of possible partitionings by running the Kernighan-Lin algorithm $n_\mathrm{KL}$ times. Note that the Kernighan-Lin algorithm typically finds some partitionings multiple times, so eventually the size of the pool may be significantly smaller than $n_\mathrm{KL}$. In the second step, we compute the entropy by means of Equations~(\ref{eq:entropy-rule})--(\ref{eq:entropy}) for each partitioning in the pool and select the one with the lowest entropy.

For the hierarchical partitioning, we perform exactly the same two steps for the root node and every internal node of the given tree structure: After selecting the partitioning of such a node according to the lowest entropy in the partitioning pool, the algorithm is called recursively if the child nodes are again internal nodes. Note that the tree structure of the Boolean reaction network is not inferred by the algorithm, but has to be given as an input. We usually choose a balanced binary tree in order to keep the number of degrees of freedom for each low-rank factor as similar as possible. We also emphasize that when a reaction pathway was approximated on the parent node, then it must not be considered for the computation of the entropies in the child nodes; otherwise the entropy for this pathway would be counted twice.

Based on the results of Section~\ref{sec:experiments}, we suggest to generate a partitioning pool with at least $n_{\mathrm{KL}} \approx 1000$ Kernighan-Lin runs. Depending on the initial partitioning and the topology of the reaction graph, also a smaller number of runs might be sufficient. In Table~\ref{tab:wall-time}, we show typical run times for generating the partitioning pool and for computing the information entropy for the largest example studied in Section~\ref{sec:experiments}, the $41$-dimensional apoptosis model. The computational time required is negligible compared to even running a single simulation with the DLR algorithm.

\begin{table}[!htb]
    \centering
    \begin{tabular}{r|c|c}
        \hline
        $n_\mathrm{KL}$\quad($n_\mathrm{PP}$) & \thead{run time \\ Kernighan-Lin} & \thead{run time \\ entropy}
        \tabularnewline
        \hline
        $1000$\quad($191$) & $2\,\mathrm{s}$ & $18\,\mathrm{s}$
        \tabularnewline
        $10\,000$\quad($687$) & $13\,\mathrm{s}$ & $70\,\mathrm{s}$
        \tabularnewline
        \hline
    \end{tabular}
    \caption{Run times for generating a partitioning pool of the $41$-dimensional apoptosis model (cf. Section~\ref{sec:experiments}) with $n_\mathrm{KL}$ Kernighan-Lin runs (resulting in a pool with $n_\mathrm{PP}$ distinct partitionings) and for computing the information entropy $H$ of all distinct partitionings in the pool. All computations were performed on a MacBook Pro with a $2\,\mathrm{GHz}$ Intel Core i5 Skylake (6360U) processor. For the Kernighan-Lin algorithm we used the implementation of the \texttt{NetworkX} package \cite{Hagberg_2008}.}
    \label{tab:wall-time}
\end{table}

We summarize the automatic partitioning scheme for the hierarchical partitioning in Algorithm~\ref{alg:automatic-partitioning}. It has to be called for the root node to obtain the partitioning for the specified tree structure.

\begin{algorithm}[!htb]
    \caption{\label{alg:automatic-partitioning}}
    \begin{tabular}{ll}
        \textbf{Input:} & Root or internal node $\tau$, number of Kernighan-Lin runs $n_{\mathrm{KL}}$ \\
        \textbf{Output:} & Partitioning $(\mathcal{P}^{\tau_0}, \mathcal{P}^{\tau_1})$\\
    \end{tabular}
    \begin{algorithmic}[1]
    \State Generate a pool of distinct partitionings by running the Kernighan-Lin algorithm $n_{\mathrm{KL}}$ times
    \State Compute total entropy $H$ using Equations~(\ref{eq:entropy-rule}) to (\ref{eq:entropy}) for each partitioning in the pool
    \State Select the partitioning $(\mathcal{P}^{\tau_0}, \mathcal{P}^{\tau_1})$ with the lowest entropy $H$
    \If{$\tau_0$ is an internal node}
    \State Call \Call{Algorithm~\ref{alg:automatic-partitioning}}{} for node $\tau_0$
    \EndIf
    \If{$\tau_1$ is an internal node}
    \State Call \Call{Algorithm~\ref{alg:automatic-partitioning}}{} for node $\tau_1$
    \EndIf
    \end{algorithmic}
\end{algorithm}

\section{Numerical experiments}\label{sec:experiments}
We tested our automatic partitioning scheme with three models from biochemistry. First, we studied the one-level partitioning of a small (a model of the mTOR pathway with $22$ species) and a medium sized reaction network (a model of signaling pathways in pancreatic cancer with $34$ species). In a second step, we examined our approach for hierarchical partitionings again with the medium sized network and eventually with a large reaction network consisting of $41$ species (a model of apoptosis). The theoretical memory requirements for solving the full CME and for the DLR approximation of each model are listed in Table~\ref{tab:memory-requirements}.

For the time integration of the low-rank quantities, we use the implementations of \cite{Prugger_2023} for one-level and \cite{Einkemmer_2024b} for hierarchical partitioning. The latter makes use of the low-rank framework \texttt{Ensign} \cite{Cassini_2021}. In both implementations we used the explicit Euler method to solve the evolution equations. We generated the partitioning pools with the \texttt{NetworkX} implementation of the Kernighan-Lin algorithm \cite{Hagberg_2008}. The reaction graphs have been plotted with \texttt{Graphviz} \cite{Ellson_2001}.

\begin{table}[!htb]
    \begin{subtable}[t]{0.49\textwidth}
        \caption*{\textbf{mTOR pathway (one-level)}}
        \centering
        \begin{tabular}{c|c|c}
            \hline
            & memory & compression ratio
            \tabularnewline
            \hline
            full & $34\,\mathrm{MB}$ & ---
            \tabularnewline
            $r=2$ & $66\,\mathrm{kB}$ & $1.94 \cdot 10^{-3}$
            \tabularnewline
            $r=4$ & $131\,\mathrm{kB}$ & $3.85 \cdot 10^{-3}$
            \tabularnewline
            $r=8$ & $263\,\mathrm{kB}$ & $7.73 \cdot 10^{-3}$
            \tabularnewline
            $r=16$ & $526\,\mathrm{kB}$ & $1.55 \cdot 10^{-2}$
            \tabularnewline
            \hline
        \end{tabular}
    \end{subtable}
    \hfill
    \begin{subtable}[t]{0.49\textwidth}
        \centering
        \caption*{\textbf{Apoptosis (hierarchical, two levels)}}
        \begin{tabular}{c|c|c}
            \hline
            & memory & compression ratio
            \tabularnewline
            \hline
            full & $18\,\mathrm{TB}$ & ---
            \tabularnewline
            $r=5$ & $207\,\mathrm{kB}$ & $1.15 \cdot 10^{-8}$
            \tabularnewline
            $r=10$ & $426\,\mathrm{kB}$ & $2.37 \cdot 10^{-8}$
            \tabularnewline
            $r=20$ & $950\,\mathrm{kB}$ & $5.28 \cdot 10^{-8}$
            \tabularnewline
            \hline
        \end{tabular}
    \end{subtable}
    \\[6ex]
    \begin{subtable}[t]{0.49\textwidth}
        \centering
        \caption*{\textbf{Pancreatic cancer (one-level)}}
        \begin{tabular}{c|c|c}
            \hline
            & memory & compression ratio
            \tabularnewline
            \hline
            full & $137\,\mathrm{GB}$ & ---
            \tabularnewline
            $r=5$ & $11\,\mathrm{MB}$ & $8.03 \cdot 10^{-5}$
            \tabularnewline
            $r=10$ & $21\,\mathrm{MB}$ & $1.53 \cdot 10^{-4}$
            \tabularnewline
            $r=20$ & $42\,\mathrm{MB}$ & $3.07 \cdot 10^{-4}$
            \tabularnewline
            \hline
        \end{tabular}
    \end{subtable}
    \hfill
    \begin{subtable}[t]{0.49\textwidth}
        \caption*{\textbf{Pancreatic cancer (hierarchical, two levels)}}
        \centering
        \begin{tabular}{c|c|c}
            \hline
            & memory & compression ratio
            \tabularnewline
            \hline
            full & $137\,\mathrm{GB}$ & ---
            \tabularnewline
            $r=5$ & $64\,\mathrm{kB}$ & $4.67 \cdot 10^{-7}$
            \tabularnewline
            $r=10$ & $140\,\mathrm{kB}$ & $1.02 \cdot 10^{-6}$
            \tabularnewline
            $r=20$ & $377\,\mathrm{kB}$ & $2.75 \cdot 10^{-6}$
            \tabularnewline
            \hline
        \end{tabular}
    \end{subtable}
    \caption{Theoretical memory requirements and compression ratios for storing either the full probability distribution (``full") or the DLR approximation for a given rank $r$ for the mTOR and pancreatic cancer models with one-level partitioning (left) as well as the pancreatic cancer and apoptosis examples with hierarchical partitioning using two levels and with the same rank for all nodes (right). The compression ratio was computed as the ratio between the memory requirements for the DLR approximation and the full probability distribution. The memory requirements are shown for a single instance of the distribution function.\label{tab:memory-requirements}}
\end{table}





\subsection{One-level partitioning}

\subsubsection{mTOR pathway}
\begin{figure}[!htb]
    \centering
    \includegraphics[width=0.6\textwidth]{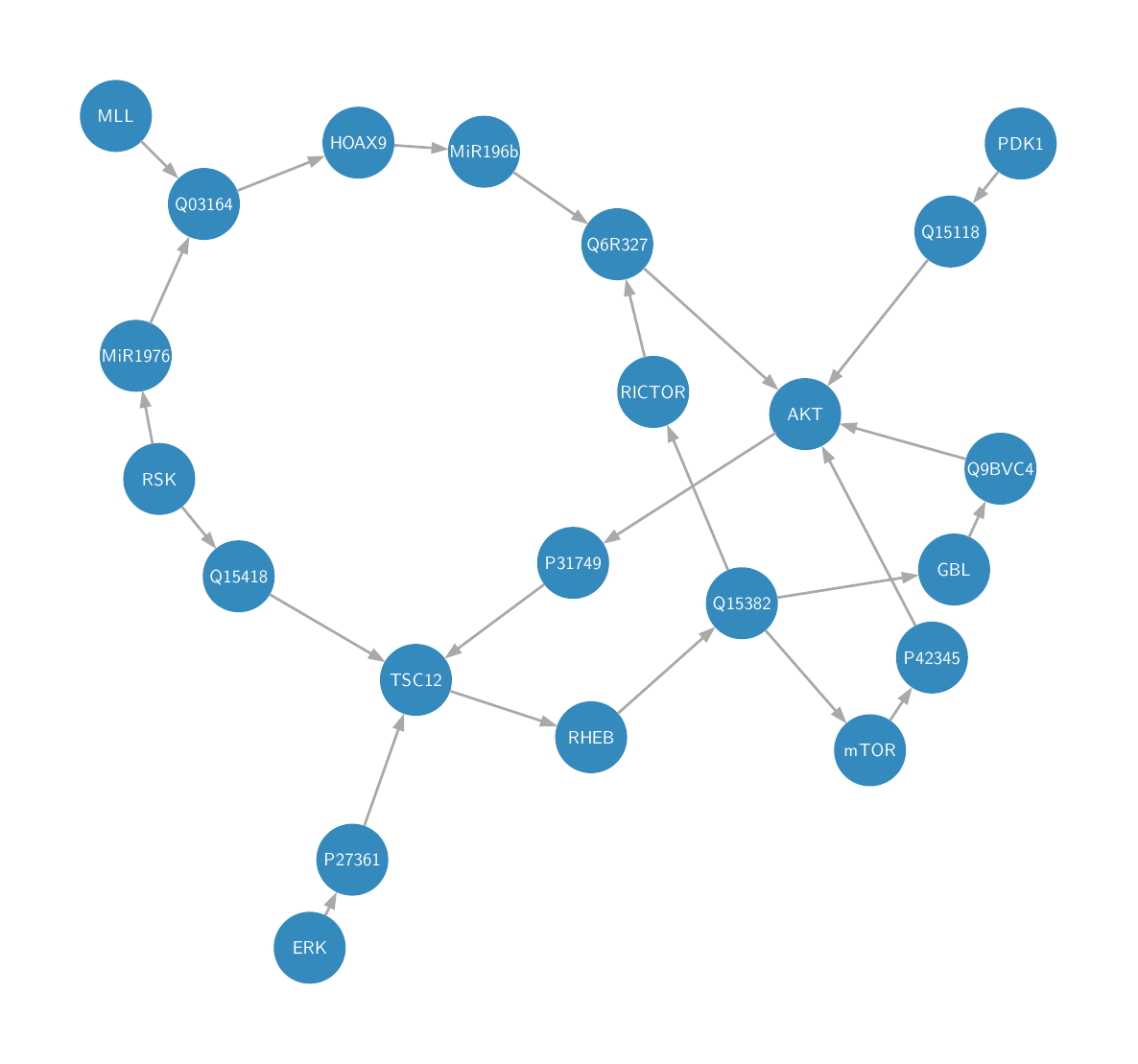}
    \caption{Reaction graph of the mTOR pathway example consisting of $22$ species.\label{fig:mTOR-graph}}
\end{figure}

In the first example, we test our partitioning scheme with a subset of the mTOR pathway suggested by \cite{Benso_2014}, which consists of $22$ species. An impairment of the mTOR pathway can lead to various diseases and targeting this pathway with drugs may be a way to treat for example leukemia or pancreatic cancer. We show the corresponding reaction graph of this example in Figure~\ref{fig:mTOR-graph}. 

Due to the relatively small system size, a reference solution could be computed by solving the master equation for Boolean reaction networks directly. Thus, for this problem we can compare the DLR simulation directly with the ground truth. We emphasize, however, that running the DLR algorithm requires significantly less memory and computational time.

As a first step, we generated a pool of one-level partitionings by running the Kernighan-Lin algorithm $100$ times. In total, $32$ different partitionings were found. In particular, three out of these $32$ partitionings (called ``PURPLE", ``BLUE" and ``GREEN") have already been studied in \cite{Prugger_2023}, where they were proposed as partitionings guided by expert knowledge. We compared these three partitionings in more detail with five partitionings which were most commonly found by the Kernighan-Lin algorithm. Furthermore, we also considered the M32 partitioning, which has the highest number of approximated reaction pathways (cuts) found by the Kernighan-Lin algorithm as a comparison. Note that the naming scheme ranks the partitionings by how often they have been found. For instance, M1 was found most often, whereas M32 was found only once (according to this naming convention, the partitionings BLUE, PURPLE and GREEN correspond to M2, M5 and M12, respectively). We integrated these nine partitions in the time interval $[0,300]$ using a time step size of $\Delta t = 10^{-1}$, with ranks $r=2$, $4$, $8$ and $16$. We chose uniform initial conditions, i.e. 
\begin{equation}\label{eq:ic-mTOR}
    P(0,\vb{x}) = \frac{1}{2^d},
\end{equation}
where $d = 22$. Note that the initial conditions are low-rank with rank $r=1$.

\begin{figure}[!htb]
    \centering
    \hspace*{0.7cm}\includegraphics[width=0.62\textwidth]{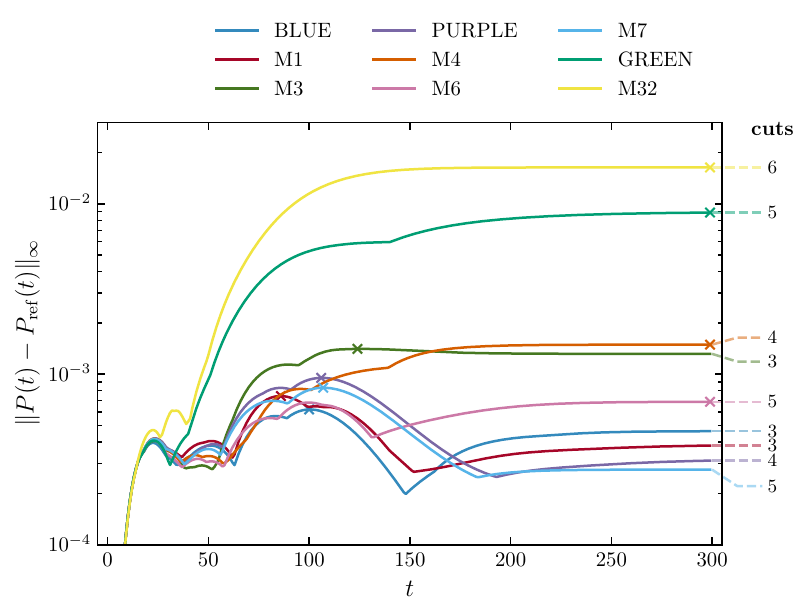}
    \caption{Time-dependent error of the DLR approximation ($r=16$) in the infinity norm for nine one-level partitionings of the mTOR pathway example. The maximum error over time for each partitioning is highlighted with a cross.\label{fig:mTOR-error-t}}
\end{figure}

In Figure~\ref{fig:mTOR-error-t}, we show the time-dependent error between the DLR approximation of the probability distribution with rank $r=16$ and the reference solution in the infinity norm for the nine chosen partitionings. According to Table~\ref{tab:memory-requirements}, the DLR approximation with $r=16$ reduces the memory requirements for storing the full probability distribution by almost two orders of magnitude. At the same time, there is a difference of almost two orders of magnitude between the least and the most accurate solution in Figure~\ref{fig:mTOR-error-t}. This clearly indicates that the choice of the partitioning highly influences the accuracy of the solution for a fixed rank. Conversely, if a solution with a given accuracy is required, for example $10^{-3}$, then for at least four partitionings in Figure~\ref{fig:mTOR-error-t} the rank has to be increased.

We also annotated the results in this figure with the number of cuts of the individual partitionings. The expectation is that the partitioning with the smallest number of cuts yields the most accurate solution and the partitioning with the most cuts indeed exhibits the worst accuracy. According to Figure~\ref{fig:mTOR-error-t}, the partitioning M32 with the largest number of cuts indeed shows the largest error, and the partitionings M1 and BLUE with only three cuts yield the best results when considering a specific time interval (around $t=150$). However, we also note that the GREEN partitioning with only three cuts is not very accurate. Moreover, when considering the accuracy of the steady states, we observe that the most accurate results are obtained from partitionings with four and five cuts (M7 and PURPLE) and not from the partitionings with the lowest number of cuts.

This substantiates the observation of Section~\ref{sec:automatic-partitioning}, that not only the number of severed reactions pathways plays a role, but also how much information is lost when cutting the pathways.

In a second experiment, we computed the entropies of all $32$ partitionings in the pool and selected the seven partitionings with the lowest entropy and also considered M32 and GREEN (partitionings with a large entropy) for further comparison. These nine partitionings were integrated in time with the same parameters as stated above.

\begin{figure}[!htb]
    \centering
    \includegraphics[width=0.57\textwidth]{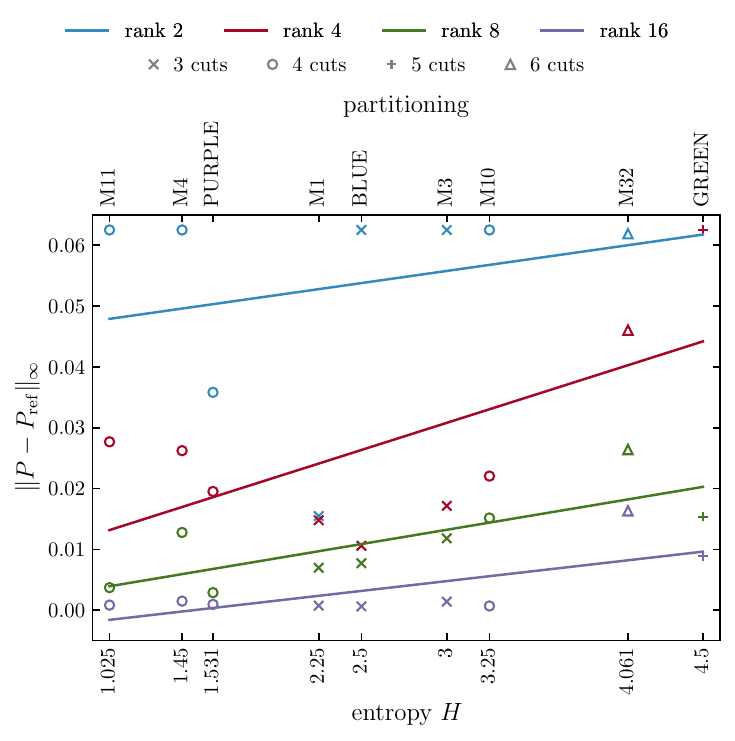} 
    \caption{Entropy-dependent error of the DLR approximation at time $t=300$ in the infinity norm for nine one-level partitionings of the mTOR pathway example. The number of cuts of the partitionings are indicated with different marker types. A linear fit is shown for solutions with the same rank.\label{fig:mTOR-correlation}}
\end{figure}

In Figure~\ref{fig:mTOR-correlation}, we plot the error for the DLR approximation of the probability distribution at time $t=300$ in the infinity norm for the nine selected partitionings with different ranks against the entropy $H$. The number of cuts of the partitionings are indicated with different marker types in the plot. We again observe that the number of cuts is not necessarily a good indicator for the resulting accuracy of the DLR approximation; for example when using $r=8$, the partitionings with the smallest number of cuts (3 cuts) do not yield the most accurate solution. In fact, in this case, there are two partitionings with 4 cuts (M11 and PURPLE) which have lower entropy than all partitionings with only 3 cuts. These two partitionings with 4 cuts outperform all partitionings with 3 cuts in terms of the error achieved.

In general, there is a correlation between the entropy and the accuracy of the DLR approximation. When choosing a small rank $r=2$ or $r=4$, all solutions exhibit a very high error, regardless of the actual partitioning, and therefore no trend between entropy and accuracy is discernible. For higher ranks such as $r=8$ the linear fit indicates that there is a correlation, and therefore, the entropy can be used as heuristic measure for selecting ``good" partitions before doing the actual simulation with the DLR approximation. We also point out that a smaller entropy does not automatically yield better accuracy for this example (compare M11 with M4). There is, however, a clear trend overall that smaller entropy partitions outperform those with comparably larger entropy (here the M32 and GREEN partitionings). Before we proceed, let us note that for the larger examples it will even be more evident that the entropy is a good measure for selecting partitionings.

\subsubsection{Pancreatic cancer}
In the second example we test out algorithm for a model of three signaling pathways which are influenced by the High-Mobility-Group-Protein B1 (HMGB1). This model was proposed in \cite{Gong_2011} and consists of $34$ species. These three pathways play an important role in many diseases, and amongst them, also in pancreatic cancer, which coins our name for this example.

We first generated a pool of one-level partitionings by running the Kernighan-Lin algorithm $10\,000$ times. As already discussed in Section \ref{sec:automatic-partitioning}, such a high number of runs is not needed in practice -- we mainly chose this number to explore exhaustively the landscape of different partitionings. Next, we computed the entropy for the distinct partititonings and selected the partitioning with the smallest entropy (\emph{best}, this partitioning would be the result of Algorithm~\ref{alg:automatic-partitioning}), the partitioning with the lowest number of cuts (which was also  the most commonly found partitioning by the Kernighan-Lin algorithm) (\emph{reasonable}) and the partitioning with the largest entropy (\emph{worst}). Additionally, we considered a partitioning which was also studied in \cite{Prugger_2023} (\emph{literature}). The four partitionings are illustrated in Figure~\ref{fig:pancreatic-graph}, the numbers of cuts and the entropy values are shown in Table~\ref{tab:pancreatic-matrix-entropy}. Note that in this figure, the first three reaction graphs show partitions (nodes which are either colored or gray, regardless of the brightness) and subpartitions (nodes of the same color and brightness). Here we only study one-level partitionings, so the subpartitions become only relevant in the next section.

In Figure~\ref{fig:pancreatic-graph} we can see that the individual partitions of the three partitionings \emph{best}, \emph{reasonable}, and \emph{literature} are contiguous (with the small exception of the INK4a node in \emph{literature}), but for the \emph{worst} partitioning, the partitions are disconnected and even form small ``islands" (such as the TLR24, HMGB1 and RAGE nodes). It is to be expected that a much higher rank is required for the DLR approximation of such a partitioning to reach a given accuracy. Remember that the partitionings were generated with the Kernighan-Lin algorithm, which tries to minimize the number of cuts between partitions. Therefore, even the \emph{worst} partitioning does not really represent the actual worst case, where the partitioning is done by randomly selecting nodes from the reaction graph.

For all four partitionings, we performed the time integration in the time interval $[0, 20]$ with time step size $\Delta t = 10^{-2}$ and ranks $r = 5$, $10$ and $20$. We again chose the uniform initial conditions of Equation~(\ref{eq:ic-mTOR}), now with $d=34$. For the benchmark we computed the first moments of the probability distribution:
\begin{equation}\label{eq:moments}
    \langle x_i \rangle (t) = \sum_{\vb{x} \in \Omega} x_i P(t, \vb{x}),
\end{equation}
where $i=0,\dots,d-1$. This can be done efficiently without computing the entire probability distribution by using Equations~(\ref{eq:DLR-approximation}) and (\ref{eq:TTN-X}).

\begin{table}[H]
    \centering
    \begin{tabular}{c|c|c}
        \hline
        partitioning & cuts & entropy $H$
        \tabularnewline
        \hline
        best & $8$ & $1.269$
        \tabularnewline
        reasonable & $6$ & $2.312$
        \tabularnewline
        worst & $12$ & $6.39$
        \tabularnewline
        literature & $10$ & $2.01$
        \tabularnewline
        \hline
    \end{tabular}
    \caption{Entropy $H$ and number of severed reaction pathways (``cuts") for four one-level partitionings of the pancreatic cancer example found by the Kernighan-Lin algorithm.}
    \label{tab:pancreatic-matrix-entropy}
\end{table}

Due to the large system size, computing an exact reference solution of the full master equation was no longer feasible. Instead, we computed the DLR approximation of $P(t,\vb{x})$ with time step size $\Delta t = 10^{-2}$ using the \emph{literature} partitioning with ranks $r=50$ and $r=60$ and compared the error of the first moments with respect to the infinity norm. The two solutions differ only by $3.34 \cdot 10^{-3}$, therefore the solution with $r=60$ is sufficiently converged and was chosen as the reference solution.

\begin{figure}[!htb]
    \centering
    \includegraphics[width=0.73\textwidth]{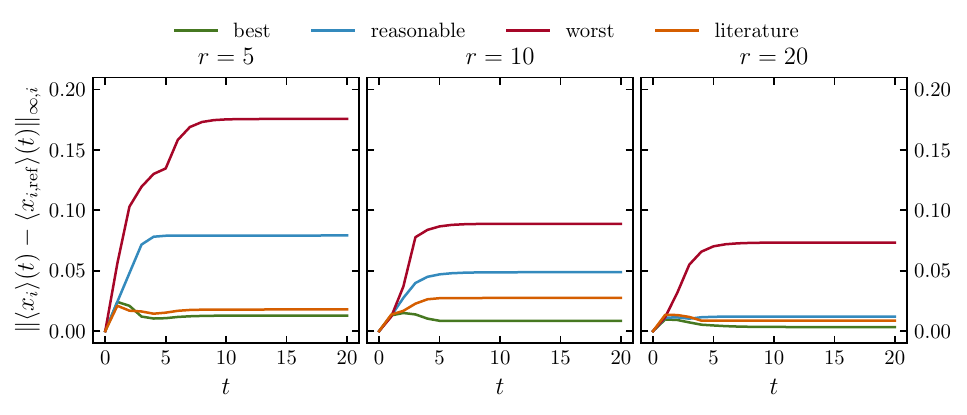}
    \caption{Time-dependent error of the DLR approximation in the infinity norm for four one-level partitionings of the pancreatic cancer example. Each subpanel shows the solutions for a different rank of the DLR approximation.\label{fig:pancreatic-matrix-moments-err}}
\end{figure}

In Figure~\ref{fig:pancreatic-matrix-moments-err}, we show the time-dependent error of the first moments for the DLR approximation in the infinity norm for the four partitionings. Each subpanel corresponds to a different rank of the DLR approximation. In this example, the entropy exactly corresponds to the observed accuracy of the solutions: the overall error for the \emph{best} partitioning is already for $r=5$ very low. At the other end of the spectrum, the \emph{worst} partitioning exhibits even for $r=20$ a relatively high error. The \emph{reasonable} partitioning (which has the lowest number of cuts) shows convergence only for $r=20$, so for this example the number of cuts is not a good measure for selecting the partitioning. Let us also point out that the partitioning chosen by Algorithm~\ref{alg:automatic-partitioning} is more accurate than the \emph{literature} partitioning chosen in \cite{Prugger_2023} by expert knowledge.

In summary, the solution of the \emph{best} partitioning proposed by our algorithm with rank $r=5$ shows already a similar accuracy as the solutions of all the other partitionings with rank $r=20$. Table~\ref{tab:memory-requirements} indicates that the memory requirements for $r=5$ and $r=20$ differ by almost a factor of $4$. Thus, selecting the partitioning with Algorithm~\ref{alg:automatic-partitioning} can significantly reduce the run time and the memory requirements of the actual simulation.

\subsection{Hierarchical partitioning}

\subsubsection{Pancreatic cancer}
In this section, we revisit the pancreatic cancer example and study hierarchical partitionings for it. To this end, we generate partitionings with Algorithm~\ref{alg:automatic-partitioning} using the same two-level, binary tree structure as shown in Figure~\ref{fig:example-tree-p1}. For each subdivision, we run the Kernighan-Lin algorithm $10\,000$ times.

The algorithm selected the partitioning with the lowest entropy (\emph{best}). Additionally, we manually selected partitionings with the lowest number of cuts (and which are at the same time also the most commonly found partitionings by the Kernighan-Lin algorithm) (\emph{reasonable}) and with the highest entropy (\emph{worst}) on every level. In Table~\ref{tab:pancreatic-hierarchical-entropy}, the number of cuts and the total entropy for each of the three partitionings is given. For hierarchical partitionings, we define the specified entropy to be the sum of the entropies of all levels. The partitionings \emph{best} and \emph{reasonable} have the same number of cuts, whereas the \emph{worst} partitioning has almost twice as many.

\begin{table}[H]
    \centering
    \begin{tabular}{c|c|c}
        \hline
        partitioning & cuts & entropy $H$
        \tabularnewline
        \hline
        best & $12$ & $3.019$
        \tabularnewline
        reasonable & $12$ & $3.749$
        \tabularnewline
        worst & $23$ & $10.859$
        \tabularnewline
        \hline
    \end{tabular}
    \caption{Entropy and number of severed reaction pathways (``cuts") for three hierarchical partitionings (using two levels) of the pancreatic cancer example.\label{tab:pancreatic-hierarchical-entropy}}
\end{table}

\begin{figure}[!htb]
    \centering
    \begin{subfigure}[b]{0.49\textwidth} 
        \centering
        \includegraphics[width=\textwidth]{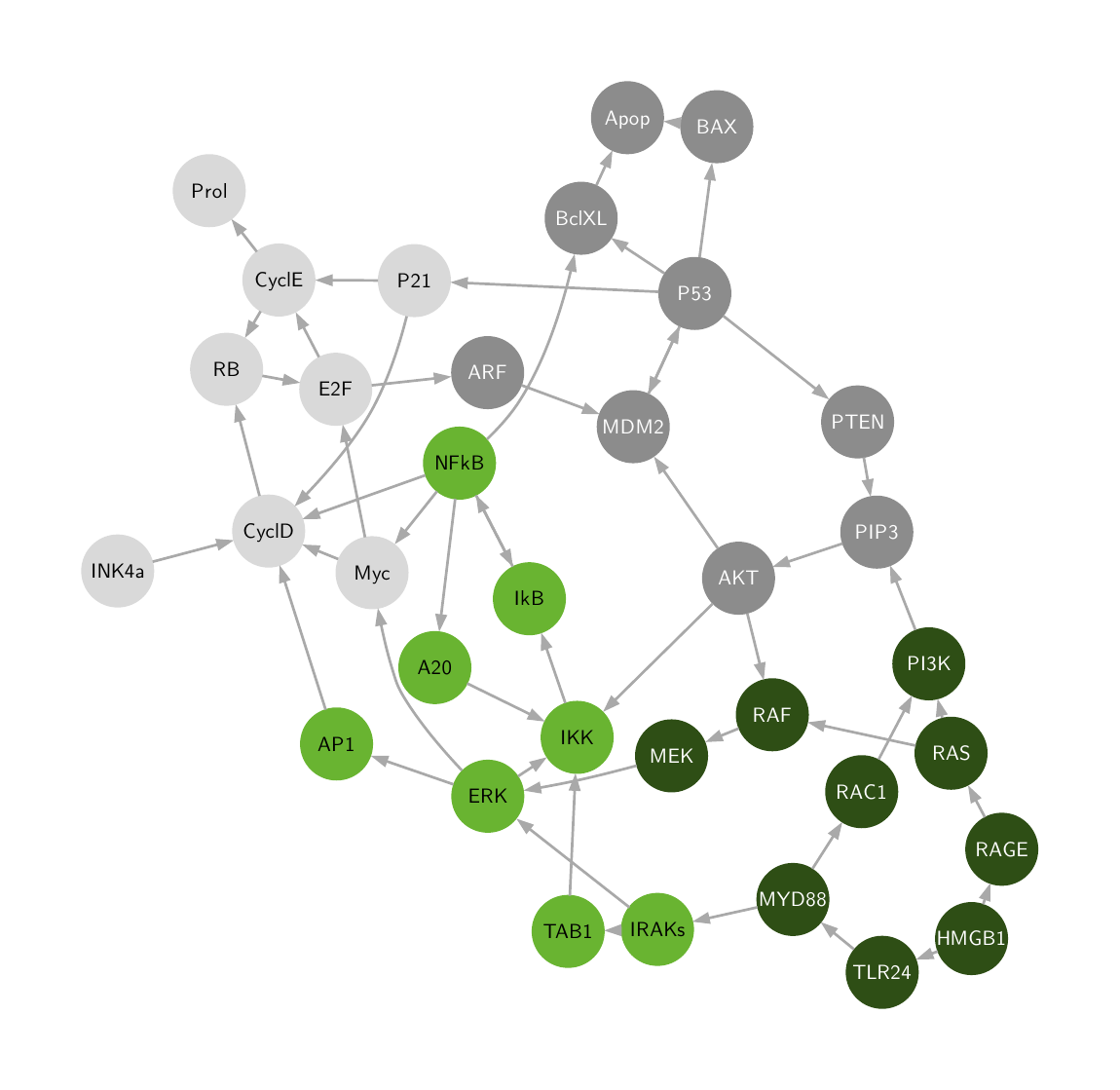}
        \caption{\label{fig:pancreatic-graph-best}}
    \end{subfigure}
    \hfill
    \begin{subfigure}[b]{0.49\textwidth}
        \centering
        \includegraphics[width=\textwidth]{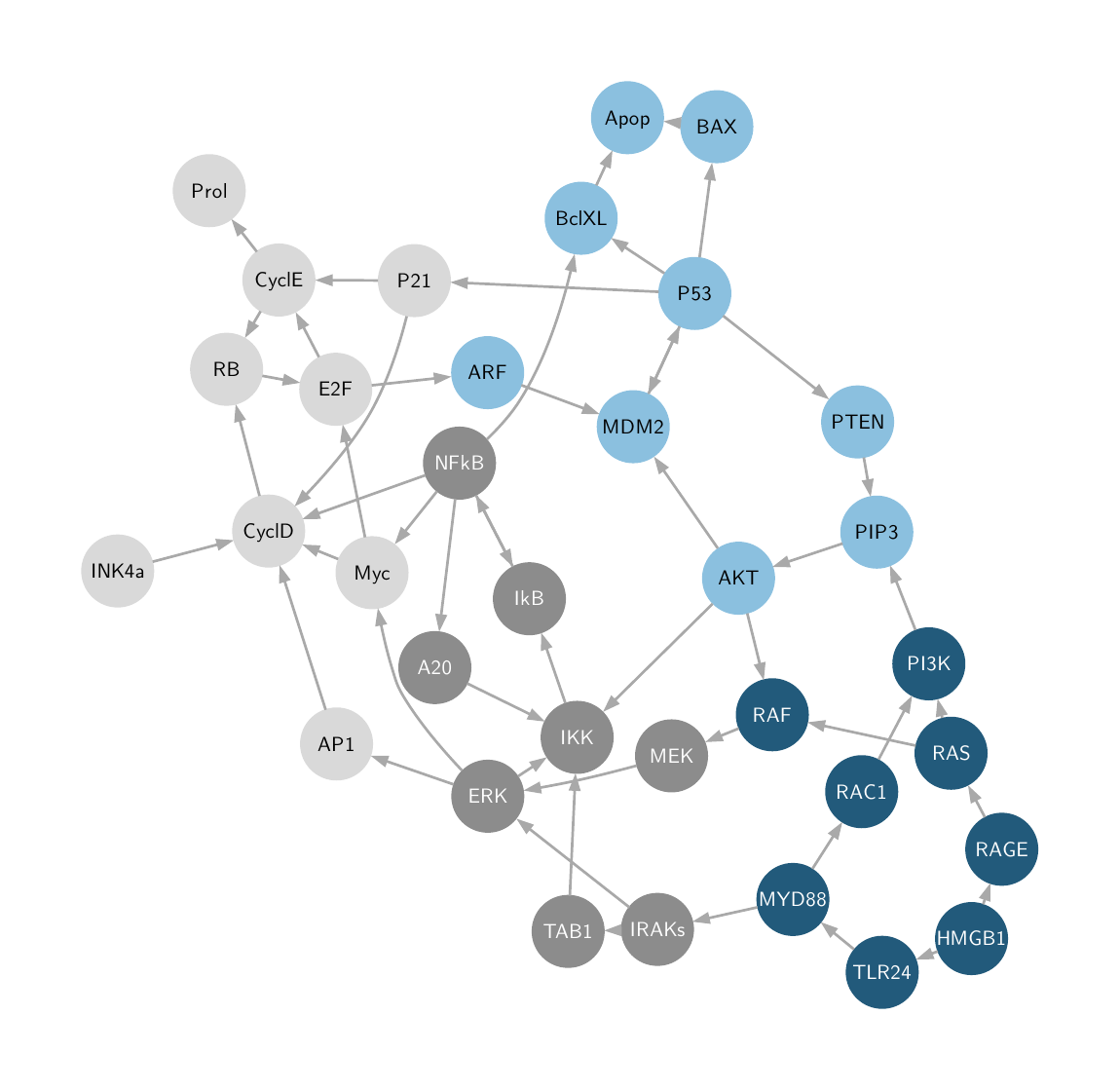}
        \caption{\label{fig:pancreatic-graph-reasonable}}
    \end{subfigure}
    \\[3ex]
    \begin{subfigure}[b]{0.49\textwidth} 
        \centering
        \includegraphics[width=\textwidth]{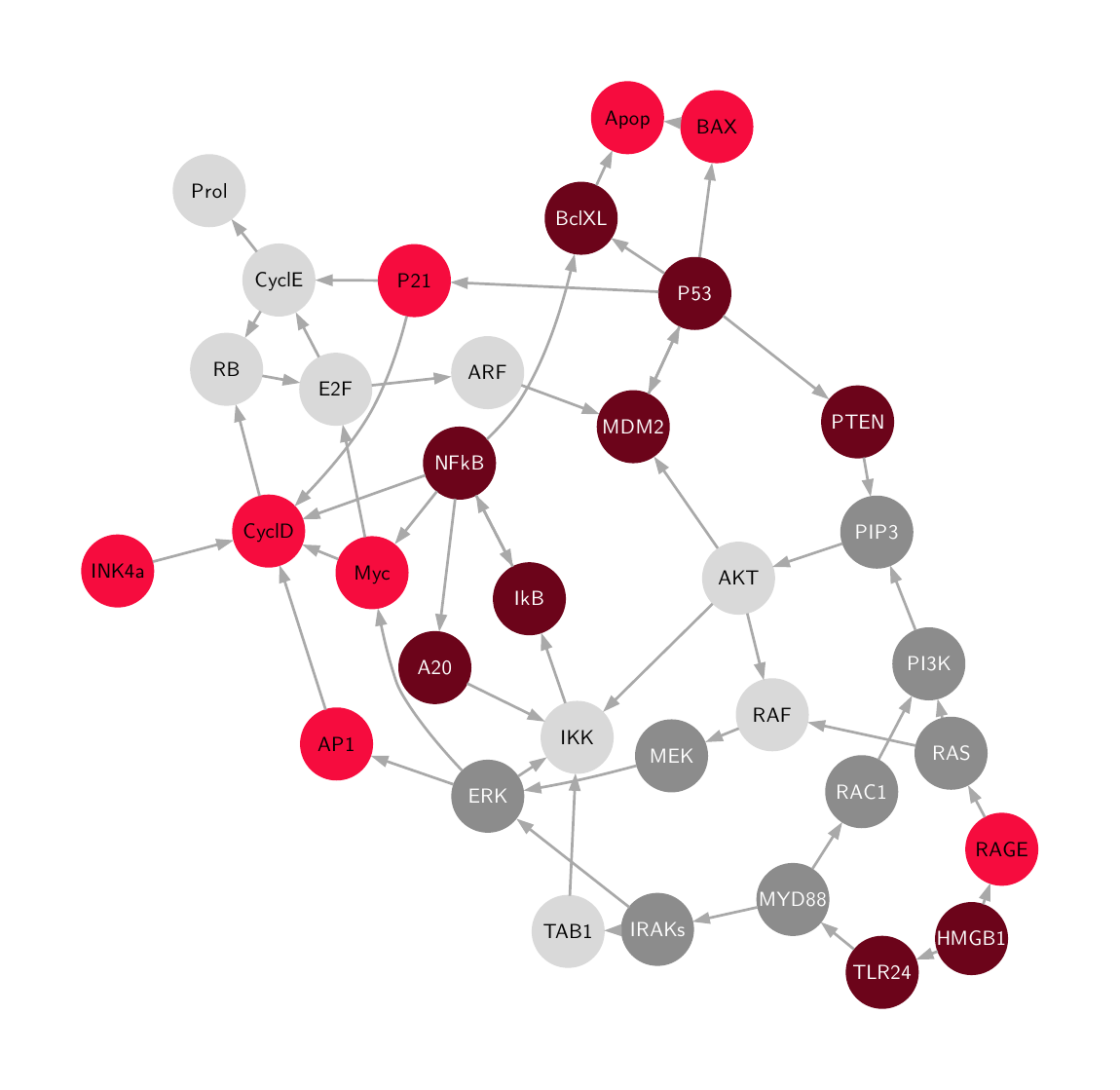} 
        \caption{\label{fig:pancreatic-graph-worst}}
    \end{subfigure}
    \hfill
    \begin{subfigure}[b]{0.49\textwidth} 
        \centering
        \includegraphics[width=\textwidth]{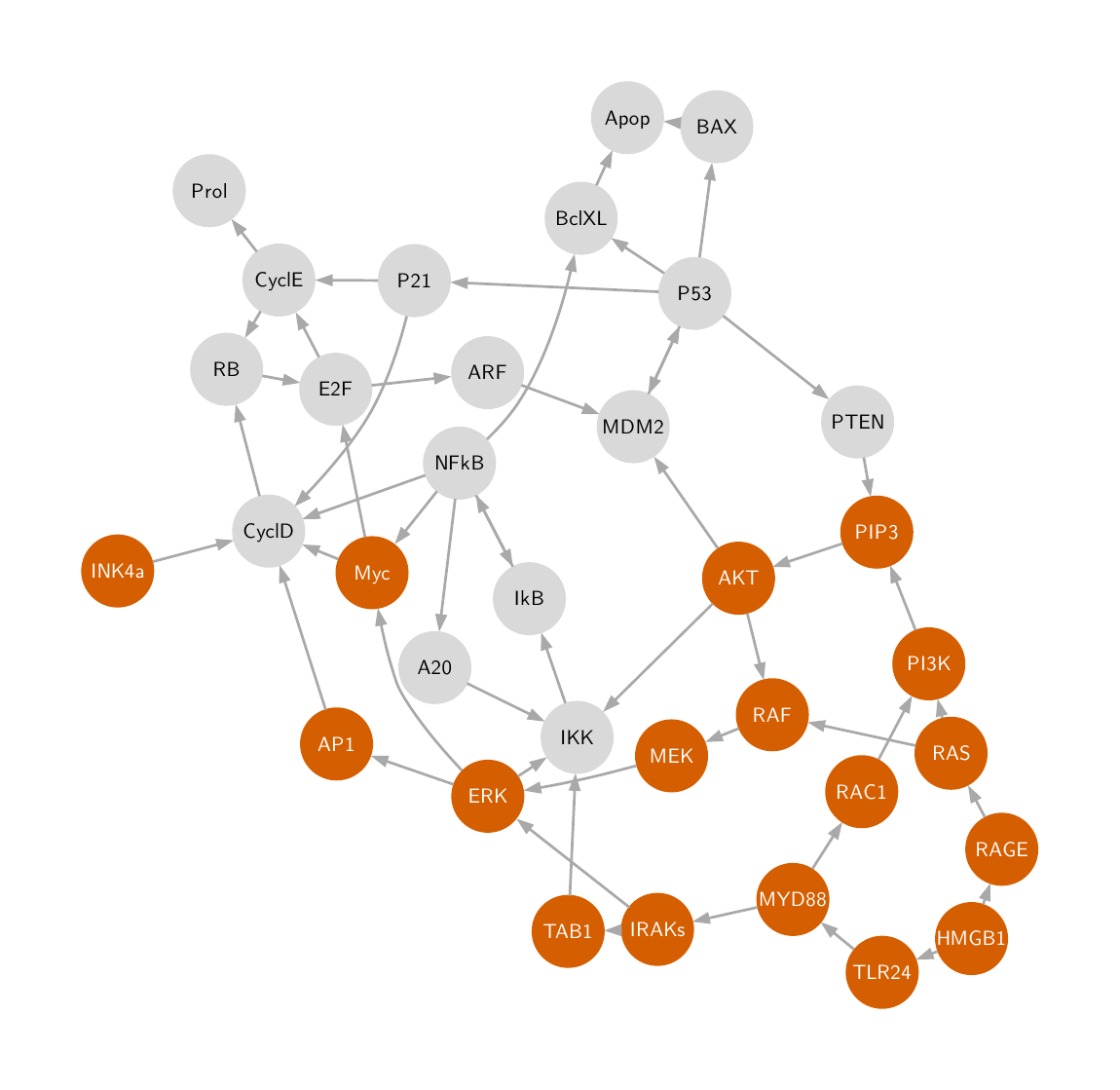} 
        \caption{\label{fig:pancreatic-graph-literature}}
    \end{subfigure}
    \caption{Reaction graphs for the pancreatic cancer example (consisting of $34$ species) with different partitionings. The colored nodes indicate to which partition and subpartition the corresponding species belong. The three graphs (a)--(c) are partitioned in a hierarchical way with the same tree structure as in Figure~\ref{fig:example-tree-p1}. The nodes are on the first level divided into a colored and a gray partition (regardless of the brightness). On the second level, the nodes are separated into ``light" and ``dark" subpartitions of the same color. (a) The partitioning with the lowest entropy (\emph{best}), (b) the partitioning with the lowest number of cuts (and which was the most commonly found by the Kernighan-Lin algorithm) (\emph{reasonable}) and (c) the partitioning with the highest entropy (\emph{worst}). (d) This reaction graph shows the one-level partitioning which was also studied in \cite{Prugger_2023} (hence we call it \emph{literature}).\label{fig:pancreatic-graph}}
\end{figure}

The reaction graphs of the three partitionings are depicted in Figure~\ref{fig:pancreatic-graph}. Note that the partitions are the same as for the one-level case, but now we also consider subpartitions (which are indicated by different brightness in the graph plots). We observe that for the \emph{best} and \emph{reasonable} partitionings also the subpartitions are contiguous, whereas for the \emph{worst} partitioning the subpartitions are much more disconnected and sometimes form small islands (e.g., the nodes for the RAGE, Apop and BAX species). Interestingly, the partitions of \emph{best} and \emph{reasonable} are almost ``orthogonal" in the sense that the subpartitions are almost the same, but the affiliation of the subpartitions to the top-level partitions is different. For each partition they have a subpartition in common (for example, the dark green subpartition for \emph{best} and the dark blue one for \emph{reasonable}), but the second subpartition belongs to the other partition (the light green subpartition in \emph{best} for instance corresponds to the dark gray subpartition in \emph{reasonable}).

For the time integration of these three partitionings we used the same parameters and initial conditions as in the previous section and also reused the reference solution of the one-level partitionings.

The time evolution of the low-rank factors, the connection tensors and the coefficient matrix was computed with ranks $r=5$, $10$ and $20$. Note that here we use the same values for the ranks of the different partitions and subpartitions, so $r^{00}=r^{01}=r$ and $r^{10}=r^{11}=r$. For comparison with the reference solution we again calculate the first moments of the probability distribution with Equation~(\ref{eq:moments}).

\begin{figure}[!htb]
    \centering
    \includegraphics[width=0.73\textwidth]{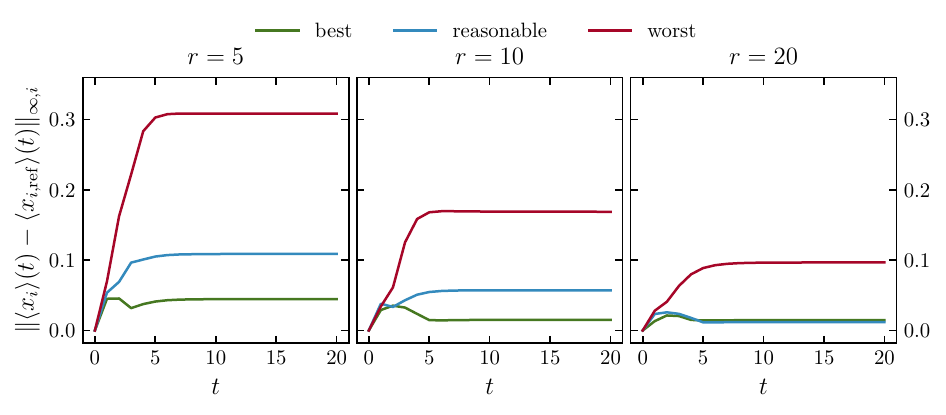}
    \caption{Time-dependent error of the DLR approximation in the infinity norm for three hierarchical (two-level) partitionings of the pancreatic cancer example. Each subpanel shows solutions for a different rank of the DLR approximation.\label{fig:pancreatic-moments-err}}
\end{figure}

The time-dependent error for the first moments of the DLR approximation in the infinity norm is displayed in Figure~\ref{fig:pancreatic-moments-err} for the three partitionings; each subpanel shows the error for a different rank.
We again find that the solution of the \emph{best} partitioning converges very quickly with increasing rank, and already for $r=5$ the error is relatively low. On the other hand, the \emph{worst} partition exhibits a very slow convergence and even for rank $r=20$ the solution is not very accurate.

Moreover, we also observe that the solutions of the hierarchical partitionings are not as accurate as the ones with one-level partitioning, but this is to be expected. The reason being that an additional approximation is introduced when we subdivide the partitions into subpartitions. Subdividing the network increases the error, but at the same time it reduces the memory and computational requirements. Table~\ref{tab:memory-requirements} shows that the memory requirements for the hierarchical partitioning with two levels are more than two orders of magnitude smaller than for one-level partitionings. The compression ratio can be up to $4.67 \cdot 10^{-7}$, which reduces the $137\,\mathrm{GB}$ for the full probability distribution to mere $64\,\mathrm{kB}$ for the DLR approximation with a hierarchical partitioning with two levels using rank $r=5$. For comparison, the one-level partitioning for the same model from the previous section needs $11\,\mathrm{MB}$ of memory and yields only a compression ratio of $8.03 \cdot 10^{-5}$ for rank $r=5$.

\begin{figure}[!htb]
    \centering
    \includegraphics[width=0.8\textwidth]{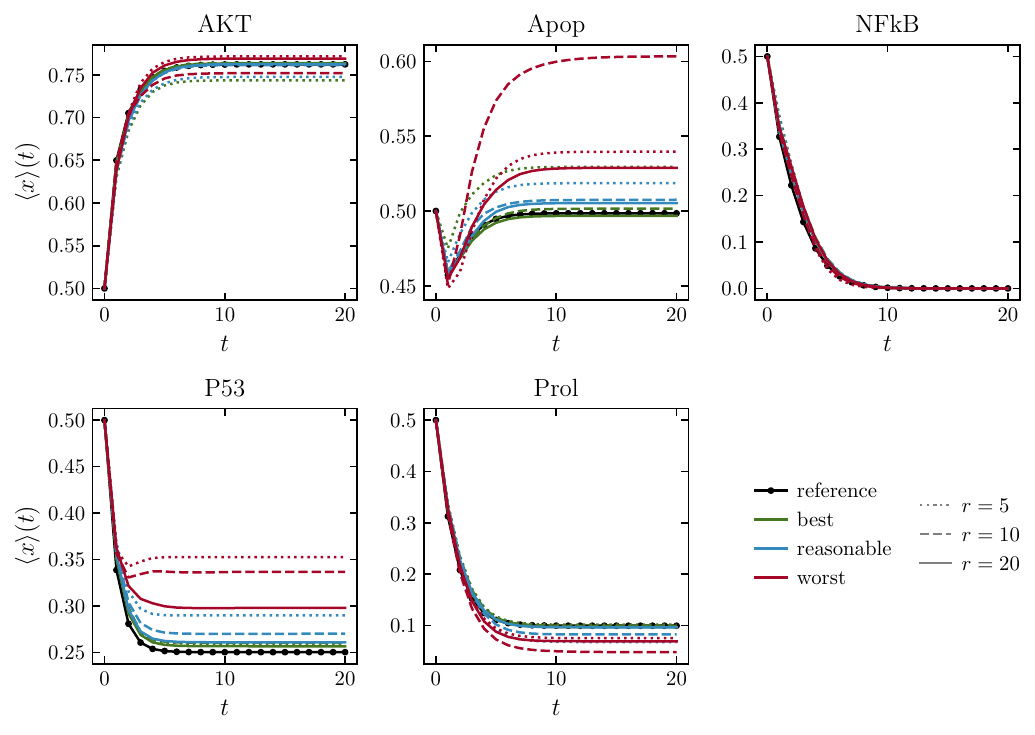} 
    \caption{Time-dependent first moments of five biochemically relevant species of the pancreatic cancer example. The plot shows the DLR approximation with three hierarchical (two-level) partitionings and the reference solution. Each color denotes a different partitioning. The DLR approximation was computed for each partitioning with three different ranks, which correspond to the different linestyles.\label{fig:pancreatic-moments-relevant-species}}
\end{figure}

In Figure~\ref{fig:pancreatic-moments-relevant-species} we show the time-dependent first moments of five biochemically relevant species for the DLR approximation and for the reference solution. We see that the overall convergence varies strongly between the different species: In particular, the \emph{worst} partitioning yields very inaccurate results for the Apop and P53 species. Both species are in the immediate vicinity of another subpartition, which could be the origin of this behavior. Apop for instance depends on BclXL, which lies in the same partition, but in a different subpartition. On the other hand, NFkB and AKT also lie on the border of two subpartitions of \emph{worse} and the convergence for these two species is significantly better. It seems that the reactions for the latter two species are not as important (i.e., convey as much information) as the ones for Apop and P53.

\subsubsection{Apoptosis}
As our last and largest example, we study the apoptosis model of \cite{Mai_2009} consisting of $41$ species proposed. Apoptosis is a cell-inherent, controlled self-destruction mechanism. An insufficient amount of apoptosis leads to uncontrolled cell proliferation and may eventually cause cancer.

We proceeded in a similar way as for the pancreatic cancer example: First, we generated three hierarchical, two-level partitionings with the tree structure shown in Figure~\ref{fig:example-tree-p1}. The first partitioning was generated with Algorithm~\ref{alg:automatic-partitioning} and exhibits the lowest entropy values (\emph{best}), the second partitioning under consideration was the one with the lowest number of cuts (and which was at the same time the most often found partitioning by the Kernighan-Lin algorithm) (\emph{reasonable}), and the third selected partitioning has the highest entropy values on all levels (\emph{worst}). The entropy values and number of severed reaction pathways for the three partitionings are tabulated in Table~\ref{tab:apoptosis-hierarchical-entropy}.

\begin{table}[H]
    \centering
    \begin{tabular}{c|c|c}
        \hline
        partitioning & cuts & entropy $H$
        \tabularnewline
        \hline
        best & $18$ & $4.207$
        \tabularnewline
        reasonable & $15$ & $6.374$
        \tabularnewline
        worst & $30$ & $19.616$
        \tabularnewline
        \hline
    \end{tabular}
    \caption{Entropy and number of severed reaction pathways (``cuts") for three hierarchical partitionings (using two levels) of the apoptosis example.\label{tab:apoptosis-hierarchical-entropy}}
\end{table}

\begin{figure}[!htb]
    \centering
    \begin{subfigure}[b]{0.49\textwidth} 
        \centering
        \includegraphics[width=\textwidth]{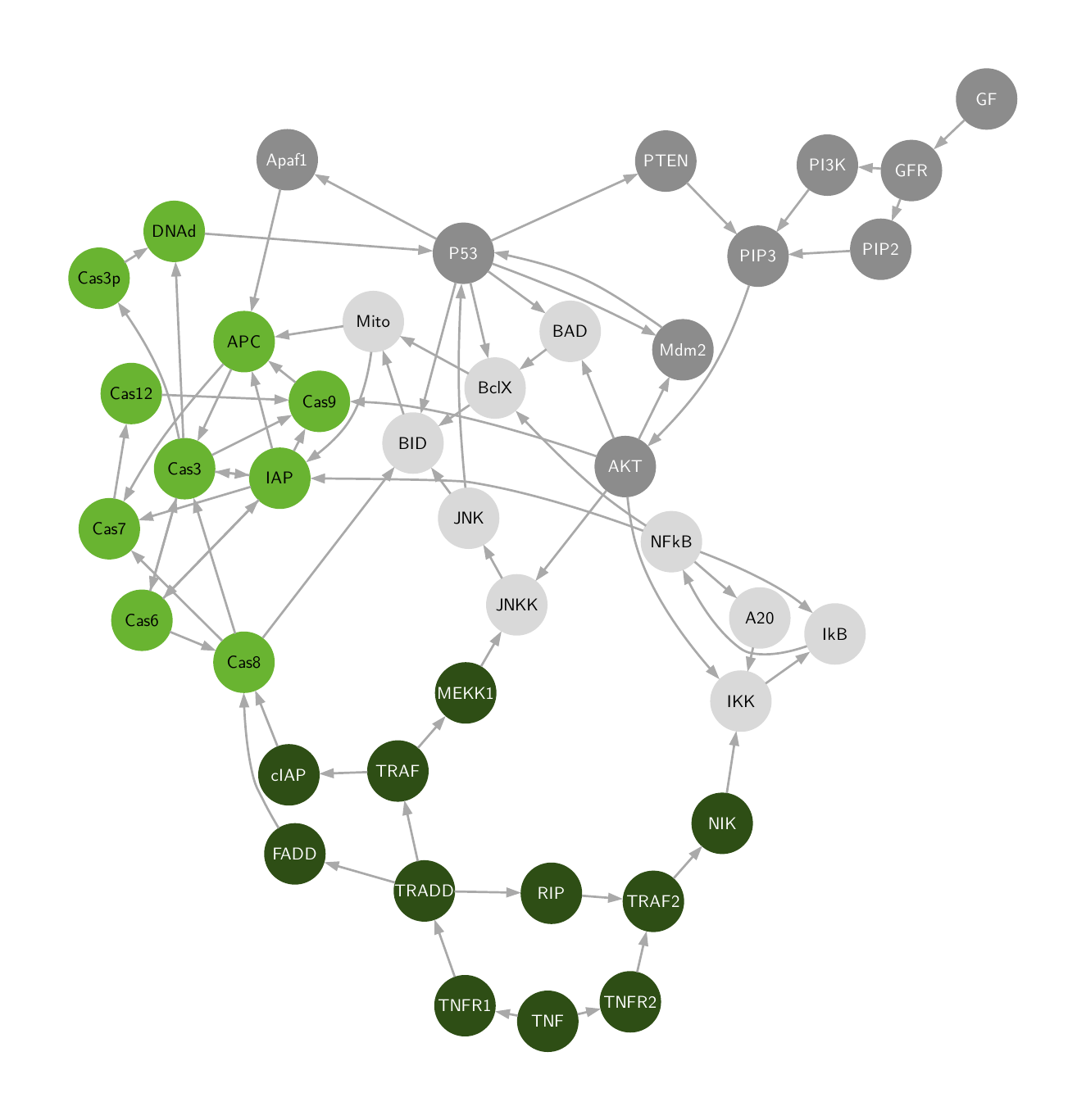}
        \caption{\label{fig:apoptosis-graph-best}}
    \end{subfigure}
    \hfill
    \begin{subfigure}[b]{0.49\textwidth}
        \centering
        \includegraphics[width=\textwidth]{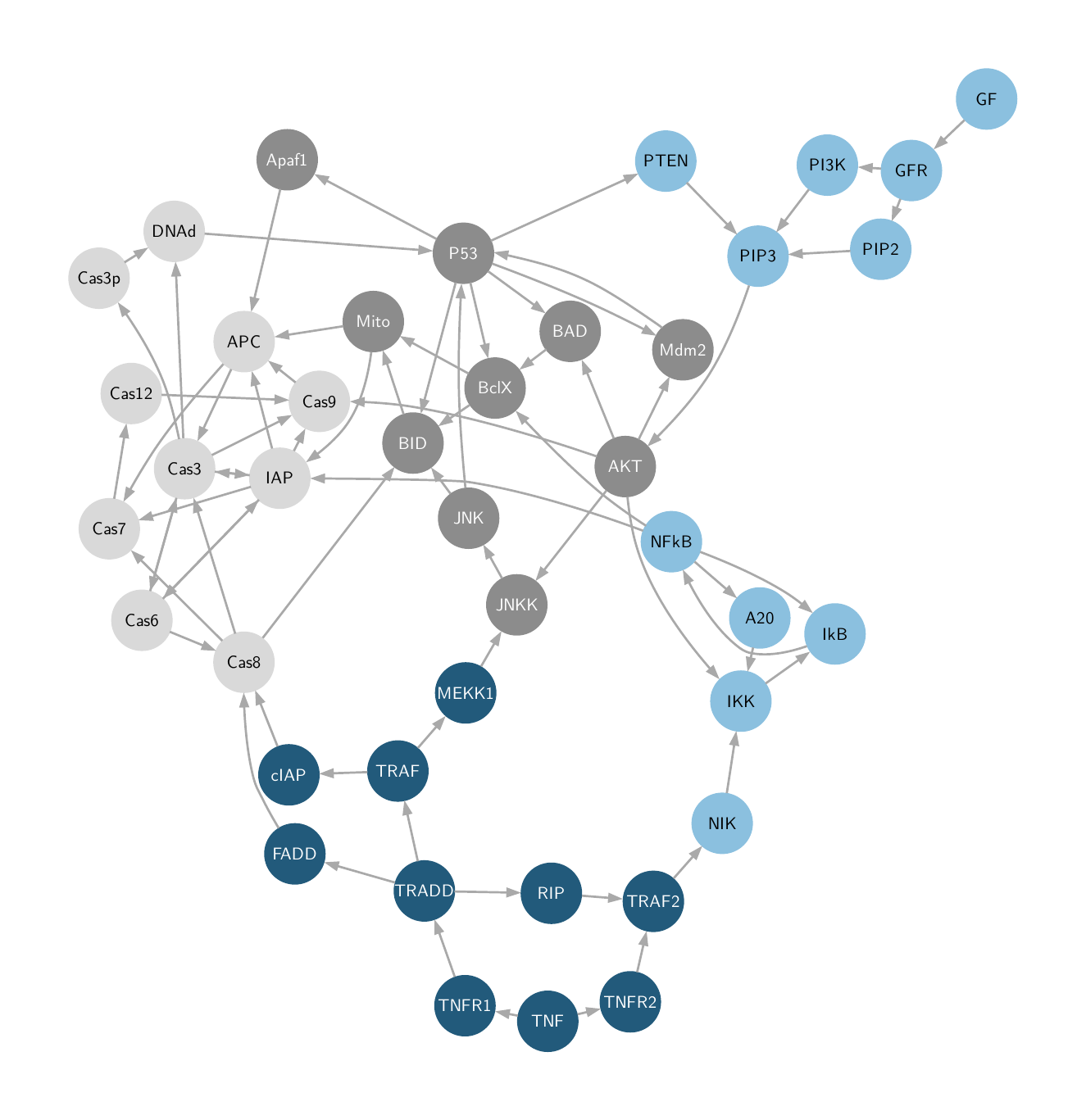}
        \caption{\label{fig:apoptosis-graph-reasonable}}
    \end{subfigure}
    \\[3ex]
    \begin{subfigure}[b]{0.49\textwidth} 
        \centering
        \includegraphics[width=\textwidth]{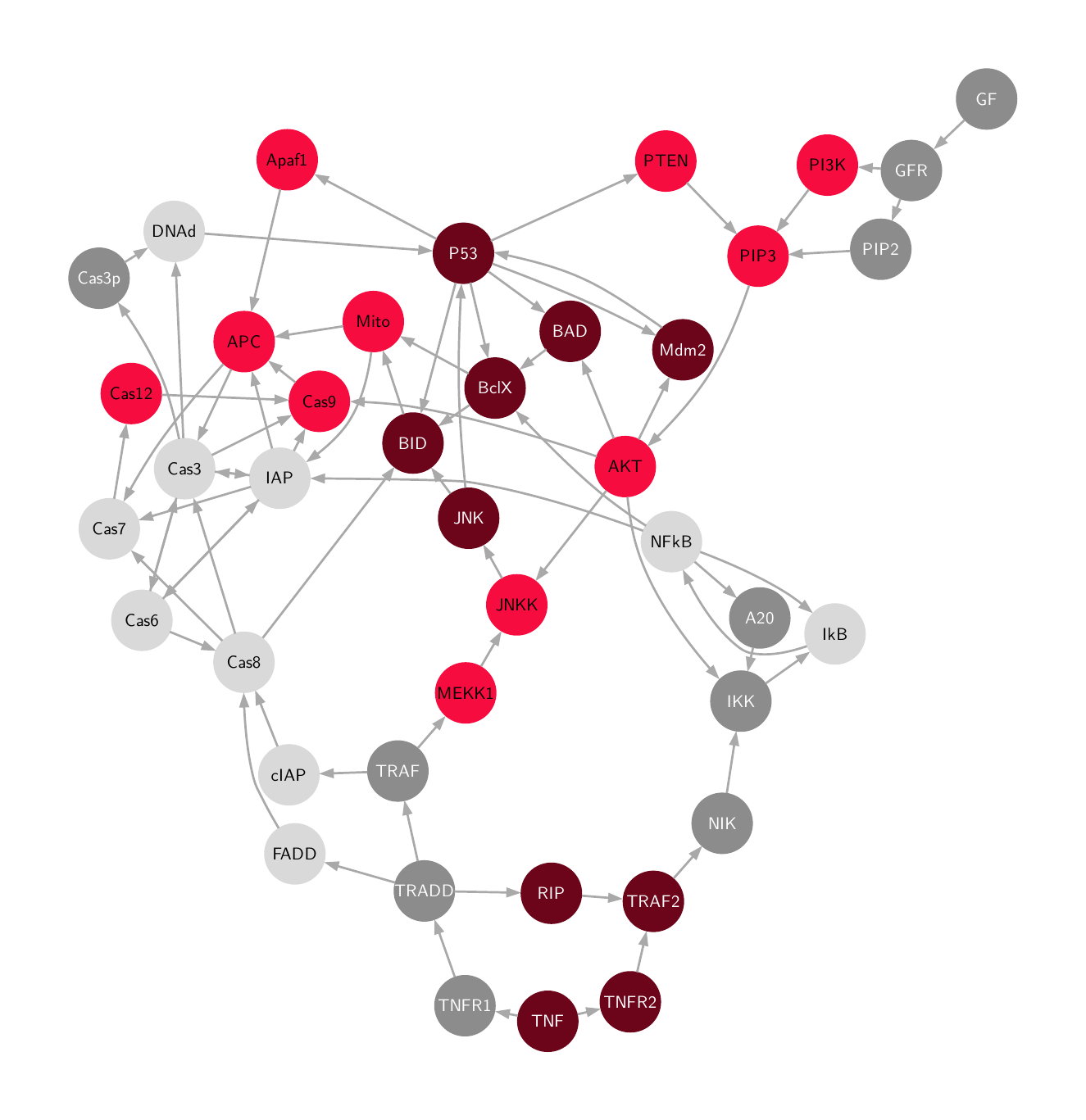} 
        \caption{\label{fig:apoptosis-graph-worst}}
    \end{subfigure}
    \caption{Reaction graphs for the apoptosis example (consisting of $41$ species) with different partitionings. The colored nodes indicate to which partition the corresponding species belong. A hierarchical partitioning scheme with the same tree structure as in Figure~\ref{fig:example-tree-p1} was used. The nodes are on the first level divided into a colored and a gray partition (regardless of the brightness). On the second level, the nodes are separated into ``light" and ``dark" subpartitions of the same color. (a) The partitioning with the lowest entropy (\emph{best}), (b) the partitioning with the lowest number of cuts (and which was the most commonly found by the Kernighan-Lin algorithm) (\emph{reasonable}) and (c) the partitioning with the highest entropy (\emph{worst}).\label{fig:apoptosis-graph}}
\end{figure}

In Figure~\ref{fig:apoptosis-graph}, the reaction graphs of the three partitionings are shown. As for the pancreatic cancer example, we again find that the \emph{worst} partitioning is very fragmented. Interestingly, we find that the dark and light green subpartitions of \emph{best} are congruent with the dark blue and the light gray subpartitions of \emph{reasonable}. Note however, that in the reaction graph for the \emph{reasonable} partitioning the light blue subpartition is separated into two islands. For the \emph{best} partitioning, all subpartitions are contiguous. We therefore expect that the \emph{best} partitioning will show better convergence than the \emph{reasonable} partitioning.

We integrate the DLR approximation for the three partitionings in the interval $[0, 20]$ with time step size $\Delta t=10^{-2}$ and ranks $r=5$, $10$ and $20$. We assume that TNF and GF are initially activated and DNAd is deactivated, and use uniform initial conditions for the remaining species:
\begin{equation*}
    P(0, \vb{x}) = \frac{1}{2^{38}} \delta_{x_\mathrm{TNF}, 1} \delta_{x_\mathrm{GF}, 1} \delta_{x_\mathrm{DNAd}, 0},
\end{equation*}
where $\delta_{x_i,x_j}$ denotes the Kronecker delta. For our benchmark, we again compute the first moments of the probability distribution with Equation~(\ref{eq:moments}).

Since computing an exact reference solution of the full master equation was again not possible due to the large system size, we computed the DLR approximation of $P(t,\vb{x})$ with time step size $\Delta t = 10^{-2}$ using a randomly selected partitioning (with a relatively low information entropy of $H=4.849$) with two levels for ranks $r=50$ and $r=60$. Comparison of the error of the first moments with respect to the infinity norm reveals that the two solutions only differ by $4.77 \cdot 10^{-3}$. Therefore, the DLR approximation with $r=60$ is sufficiently converged and was selected as the reference solution.

\begin{figure}[!htb]
    \centering
    \includegraphics[width=0.68\textwidth]{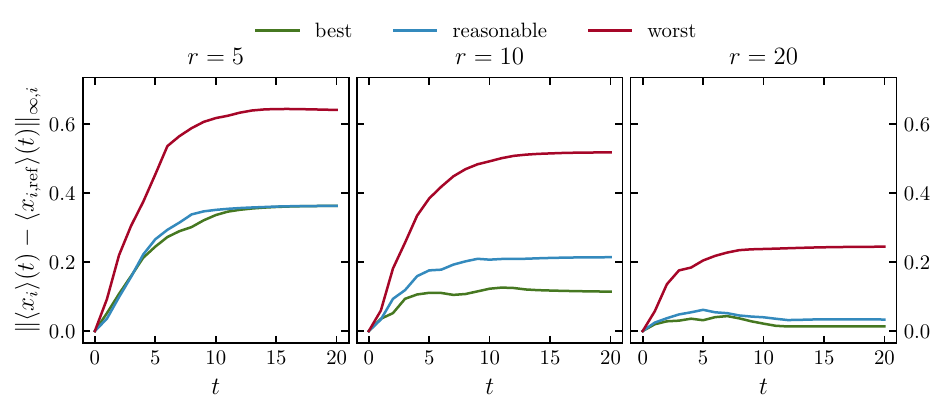} 
    \caption{Time-dependent error of the DLR approximation in the infinity norm for for three hierarchical (two-level) partitionings of the apoptosis example. Each subpanel shows solutions for a different rank of the DLR approximation.\label{fig:apoptosis-moments-err}}
\end{figure}

In Figure~\ref{fig:apoptosis-moments-err} we show the time-dependent error for the first moments of the DLR approximation in the infinity norm for the three partitionings. Similar to the previous example, we see that the \emph{worst} partitioning has the largest error for all ranks. The \emph{reasonable} and \emph{best} partitionings have a similar accuracy for $r=5$, but for higher ranks, the \emph{best} partitioning outperforms the \emph{reasonable} partitioning. Recall that the \emph{reasonable} partitioning has the lowest number of cuts, we therefore conclude that also for this example the information entropy is the better measure for selecting the partitioning.

\begin{figure}[!htb]
    \centering
    \includegraphics[width=0.78\textwidth]{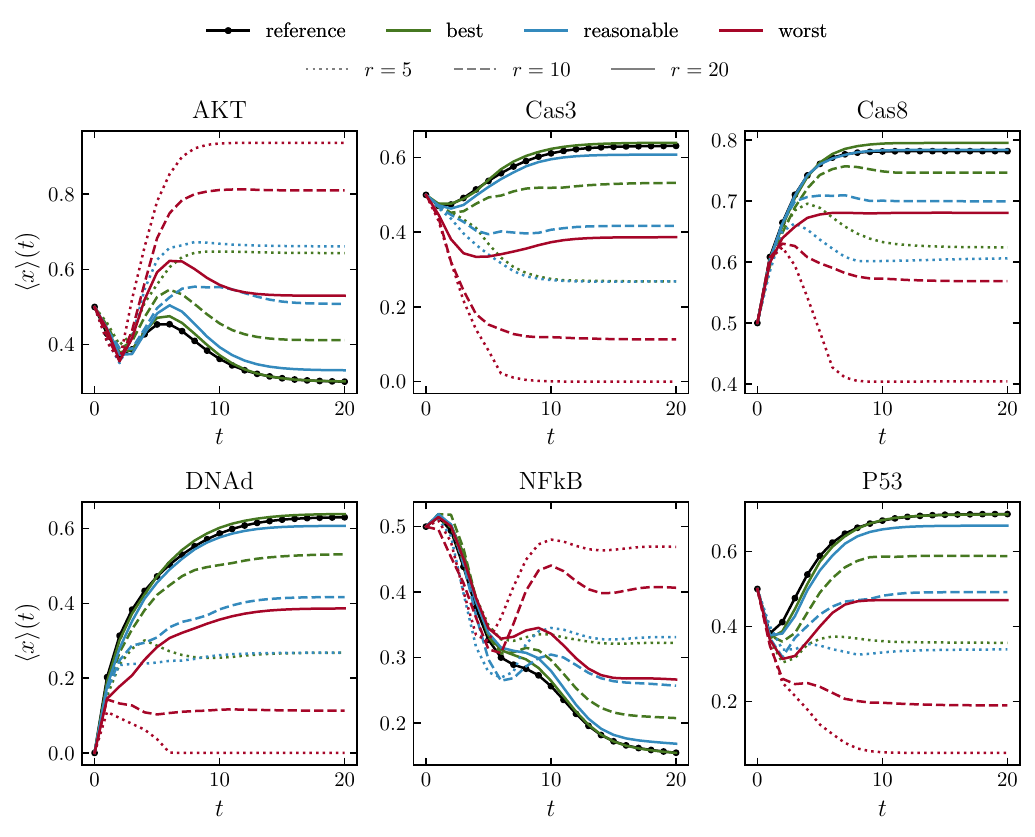} 
    \caption{Time-dependent first moments of six biochemically relevant species of the apoptosis example. The plot shows the DLR approximation with three hierarchical (two-level) partitionings and the reference solution. Each color denotes a different partitioning. The DLR approximation was computed for each partitioning with three different ranks, which correspond to the different linestyles.\label{fig:apoptosis-moments-relevant-species}}
\end{figure}

Figure~\ref{fig:apoptosis-moments-relevant-species} shows the time-dependent first moments for six biochemically relevant species for the DLR approximation and the reference solution. In contrast to the pancreatic cancer example, the overall convergence is relatively homogeneous over all species. The \emph{best} partition with $r=20$ yields the best approximation for the first moments of all species except one -- only for Cas8 the \emph{reasonable} partitioning yields a better result. Curiously, the convergence of the \emph{worst} partition is for some species very erratic. For instance for Cas8, the solution with rank $r=5$ captures the behavior of the reference solution very well, but when increasing the rank to $r=10$, the solution for the first moments shows a higher error (this is also observed for the \emph{literature} partitioning for ranks $r=5$ and $r=10$). Note that this effect of an increasing error for a higher rank is observed for the moments due to the accumulation of approximation errors; the error for the DLR approximation of the probability distribution $P(t,\vb{x})$ does decrease when a higher rank is used.

\section{Conclusion and outlook}\label{sec:conclusion}
In this work, we proposed an heuristic algorithm for the partitioning of Boolean reaction networks. The algorithm aims to find a partitioning that minimizes the rank of the DLR approximation for a given accuracy and consists of two steps. The first step employs the Kernighan-Lin algorithm, which partitions the reaction graph and minimizes the number of approximated reaction pathways. A pool of partitionings is generated by running the Kernighan-Lin algorithm several times. In the second step of the algorithm, the partitioning with the lowest information entropy is selected. The numerical experiments indicate that the partitioning selected in this way yields accuracy superior to both partitionings chosen by human experts and those found by simply minimizing the number of cuts. The Kernighan-Lin algorithm and the calculation of the entropy are computationally cheap compared to the time evolution of the DLR approximation and therefore can be implemented as a fast preprocessing step in the actual simulation.

The total entropy $H$ could also be extended to more general types of reaction networks. For Boolean reaction networks we assumed that each species has a single Boolean rule, and the probability for all Boolean rules is the same. For the more general case of probabilistic Boolean reaction networks \cite{Shmulevich_2002}, where species may have more than one rule and rules are allowed to have different probabilities, the extension is straightforward: The rule-specific entropies $\tilde{h}_i(\vb{x}_i)$ would have to be weighted with the specific probabilities and an additional sum in the definition of $H$ would have to be introduced which has to run over all rules associated with species $S_i$.

It might also be possible to generalize the information entropy for general chemical reaction networks, which are described by the kinetic chemical master equation (CME). However, such an extension is not straightforward for a number of reasons. First, the two Boolean activation states have to be replaced in this model by population numbers, which can assume any integer value. Secondly, for chemical reaction networks there are in general multiple products per reaction possible, whereas for Boolean networks only a single species is activated or deactivated per rule. Lastly, the reaction rates of the Boolean networks are replaced in the more general model by so-called propensity functions depending on the population numbers, and they would also have to be considered for the calculation of the information entropy.

An alternative approach to the information entropy could be the computation of the ranks of the individual transition probabilities. This could be done efficiently, since the transition probabilities often show a Kronecker structure as in the example in Section~\ref{sec:boolean-reaction-system}. The resulting ranks could then be used as a measure for selecting the best partitioning from the pool generated by the Kernighan-Lin algorithm. In the context of quantum spin systems something similar has been done in \cite{Ceruti_2024}.

Finally, future work could address how to incorporate entropy information in the weights of the reaction graph. Thus the Kernighan-Lin algorithm would directly minimize a functional depending on the number of cuts and the entropy. With this approach the entropy would have to be calculated only once before generating the partitioning pool.

\section*{Acknowledgments}
The authors thank Luk Burchard for valuable discussions. This work has been supported, in part, by the Tiroler Nachwuchsforschungsf\"{o}rderung (TNF) under grant F.47947/5-2023.


\printbibliography

\end{document}